\documentclass[sn-mathphys]{sn-jnl}
\jyear{2023}
\usepackage{subcaption}
\usepackage[english]{minitoc}
\usepackage{multirow}
\usepackage{float}
\restylefloat{table}
\theoremstyle{thmstyleone}
\newtheorem{theorem}{Theorem}[section]
\newtheorem{proposition}[theorem]{Proposition}
\theoremstyle{thmstyletwo}

\newtheorem{remark}{Remark}
\newtheorem{lemma}{Lemma}
\theoremstyle{thmstylethree}

\begin{document}

\author*[1]{\fnm{Martial} \sur{Longla}}\email{mlongla@olemiss.edu} 
\author[1,2]{\fnm{Mous-Abou} \sur{ Hamadou}}\email{mhamadou@go.olemiss.edu}

\affil[1]{\orgdiv{Department of mathematics}, \orgname{University of Mississippi}, \orgaddress{\street{University Ave}, \city{University}, \postcode{38677}, \state{Mississippi}, \country{US}}\\ {ORCID: 0000-0002-9432-0640} }
\affil[2]{\orgdiv{Department of mathematics}, \orgname{University of Maroua}, \country{Cameroon}}

\title{Estimation problems for some perturbations of the independence copula. }

\abstract{This work provides a study of parameter estimators based on functions of Markov chains generated by some perturbations of the independence copula. We provide asymptotic distributions of maximum likelihood estimators and confidence intervals for copula parameters of several families of copulas introduced in Longla (2023). Another set of moment-like estimators is proposed along with a multivariate central limit theorem, that provides their asymptotic distributions. We investigate the particular case of Markov chains generated by sine copulas, sine-cosine copulas and the extended Farlie-Gumbel-Morgenstern copula family. Some tests of independence are proposed. A simulation study is provided for the three copula families of interest. This simulation proposes a comparative study of the two introduced estimators and the robust estimator of Longla and Peligrad (2021), showing advantages of the proposed work.}


\maketitle

\section{Introduction}
Over the past decades, many authors have worked on building copulas with various properties such as prescribed diagonal section, known association properties, perturbation and mixing. Some of these constructions can be found in Nelsen (2006).  Arakelian and Karlis (2014) studied mixtures of copulas, that were investigated for mixing by Longla (2015). Longla (2014) constructed copulas based on prescribed $\rho$-mixing properties obtained from extending results of Beare (2010) on mixing for copula-based Markov chains. Longla and al. (2022, 2022b) constructed copulas in general based on perturbations. In these two works, two perturbation methods were considered. Our current work is based on perturbations of the independence copula. These perturbations are of the form $C_D(u,v)=C(u,v)+D(u,v)$, where $C(u,v)$ is the copula that is undergoing the perturbation and $C_D(u,v)$ the perturbed copula. This method of perturbation was presented in Komornik et al (2017), where authors mention that it allows in the case of the independence copula to introduce some level of dependence. Several other authors have considered extensions of existing copula families via other methods, that are not the focus of this work, see for instance, Morillas (2005), Aas et al (2009), Klement et al (2017), Chesneau (2021), among others. 

We should mention, that trigonometric copulas that were investigated in Chesneau (2021), are very different from the perturbations that we have in this paper. This can be seen on the level of the dependence structure and the functional form of the copulas. This paper is concerned with the important question of estimation, confidence intervals and testing for copula-based Markov chains, that are generated by copula families of Longla (2023). These families were introduced, but not studied for the mentioned questions. Their mixing properties were examined, making it possible to explore central limit theorems. We base our study of MLE on the work of Billingsley (1961) for stationary Markov chains and their asymptotic distributions. Billingsley (1961) has been used in the setup of Markov chains by several authors over the years. Recently, it has been used by Hossain and Anam (2012) to study modelling of rainfalls in Bangladesh via Markov chains and a logistic regression model.  Zhang et al. (2023) used it to study asymptotic normality of estimators in an AR(1) Binomial model. Kuljus and Ranneby (2021) considered maximum spacing as approximation of relative entropy for semi-markov processes and continuous time Markov chains. A good resource for the use of copulas in estimation problems is Chen and Fan (2006).

Our paper is dedicated to the analysis of three copula families. This analysis relates mixing coefficients defined by Bradley (2007), in the context of copulas, to the study of Markov chains. Copulas are functions that naturally represent the effect of dependence in multivariate distributions. This fact is due to Sklar's theorem (see Sklar (1959)). Copulas are obtained by scalling out the effect of marginal distributions from joint distributions. In the case of continuous marginal distributions, most mixing coefficients of Markov chains are not affected by the marginal distributions (see Longla (2015)). Combining these facts with the special nature of the copula families of interest of this paper, Longla (2023) showed that these copulas generate $\psi$-mixing. It is worth mentioning that copula coefficients are eigen-values of the copula operators. They are associated to the eigen-functions $\varphi_k(x)$ present in the formula of the copula of Longla (2023). 

  The first family that we study has a density based on orthogonal cosine functions. The second is based  on  orthogonal sine and cosine functions, and the third family is based on Legendre polynomials. Orthogonality here is in the sense of $L^2(0,1)$. This last family is an extension of the Farlie-Gumbel-Morgenstern (FGM) copula family. It introduces a second dependence coefficient $\lambda_2$. The FGM has been extensively studied in the literature with several extensions (see Hurlimann (2017) or Ebaid and al. (2022)). 

In general, a multivariate central limit theorem is provided in our paper for parameter estimators of copulas of Longla (2023). This central limit theorem is based on  Kipnis and Varadhan (1986). Kipnis and Varadhan (1986) showed that for functions of a stationary reversible Markov chain with finite mean and variance, when the variance of partial sums is asymptotically linear in $n$, the central limit theorem holds. 

A general theorem of Billingsley (1961) is used to establish the limiting distribution of MLE. We show that they exist and are unique, when the true parameters are inside the support region. Billingsley proved that under some regularity conditions (which are checked for each of the copula families), the MLE is asymptotically normal with variance provided by a form of Fisher information. For practical reasons, given that the Fisher information can't be comptuted in closed form, an estimate was used. Computation time for this procedure makes the MLE less attractive than the estimator that we propose. Based on properties of the coefficients of the construction, our proposed parameter estimators have closed form asymptotic variances. 

\subsection{Definitions and comments}
We call copula in general a bivariate copula, which is defined as a bivariate joint cumulative distibution function with uniform marginals on $[0,1]$. An equivalent definition can be found in Nelsen (2006). A stationary Markov chain can be represented by a copula and its stationary or marginal distribution. The copula of a stationary Markov chain is that of any of its consecutive variables.
We say that a stationary sequence $(Y_1,\cdots, Y_n)$ with finite mean $\mu$ and variance $\sigma^2$ satisfies the central limit theorem if  $\sqrt{n}(\bar{Y}-\mu)\to N(0,\sigma_0^2)$, for some positive number $\sigma^2_0$. Kipnis and Varadhan (1986)  showed that when the variance of partial sums of functions of the Markov chain satisfies  $\sigma_f^2=\lim_{n\to\infty}n^{-1} var(S_n)$, where $S_n=\sum_{k=1}^nf(X_i)$, $Y=f(X)$, $\mu=\mathbb{E}(f(X))$, the central limit theorem holds. This theorem is combined with the Cram\'er-Wold device to prove a multivariate central limit theorem for parameter estimators in this paper. In doing so, we use the $n$-fold product of copulas to compute the limiting variance. The $n$-fold product of the copula uses partial derivatives to find a new copula. The notation $C_{,i}(u,v)$ is for the derivative of $C(u,v)$ with respect to the $i^{th}$ variable. The $n$-fold product of a copula is defined in Nelsen (2006) (or by Darsow et al. (1992)) by $C^1(u,v)=C(u,v),$  $$\text{for}\quad n>1, C^n(u,v)=C^{n-1}*C(u,v)=\int_0^1 C^{n-1}_{,2}(u,t)C_{,1}(t,v)dt.$$
Given that we are dealing with Markov chains, we recall that the fold product for $n=2$ is the joint cumulative distribution of $(X_1,X_3)$ if $(X_1,X_2,X_3)$ is a stationary Markov chain with copula $C(u,v)$ and the uniform marginal distribution (see Darsow et al. (1992)).

In general, the copula $C^{n}(u,v)$ is the joint distribution of $(X_0, X_n)$ when $X_0, \cdots, X_n$ is a stationary Markov chain generated by $C(u,v)$ and the uniform distribution. These definitions can be found in different forms in various papers, but are equivalent when the study is concerned with absolutely continuous copulas and continuous marginal distributions. This is because Sklar's theorem guarantees uniqueness of the representation of the joint distribution of a vector through its copula and marginal distributions (see Sklar(1959)). Copulas simplify the study of mixing properties of Markov chains due to the fact that for continuous random variables, the marginal distributions are scaled out in the computation; the formulas reduce to equivalent forms in terms of the uniform marginal distributions. This is one of the advantages of dealing with copulas while studying mixing.

Mixing coefficients, that we refer to, are defined in Bradley (2007), and were used in Longla (2023) to establish properties of Markov chains for the copulas that we are studying in this work. This work is not concerned by mixing in general, but rather uses the that fact, that mixing was established in an earlier work. This justifies the limitation of the background to the provided references. The mixing rate allows convergence of the variance of partial sums of functions of the Markov chains, and implies the needed ergodicity in the Theorem of Kipnis and Varadhan (1986). This theorem on its turn is used to establish central limit theorems for copula parameter estimators. In the simulation part of this work, we use a result of Longla and Peligrad (2021) to derive some estimators of copula parameter. These estimators are compared to the MLE, and a second class of estimators, that is derived based on the fact, that parameters are eigen-values of the operator induced by the copula of the Markov chain on $L^2(0,1)$.
\subsection{Structure of the paper and main results}
We propose several estimators of copula parameters for the copulas created in Longla (2023). Large sample study of these estimators is provided in general via derivation of a multivariate central limit theorem. These results are applied to 3 concrete families of copulas. We have conducted a comparative simulation study in which 4 different estimators are used. For each of the examples, we perform maximum likelihood estimation and show that the $R$ codes for our proposed estimators take less time to run than the MLE and can be preferable in large samples.
 
The rest of the paper is organized as follows. In section 2 we provide the copula  families of interest in this work along with some graphs of copula densities and level curves. In section 3 we present parameter estimation procedures for Markov chains generated by the copula families proposed by Longla (2023). This section includes our proposed estimators, central limit theorems for estimators and $\chi^2$ tests of independence of observations. A multivariate central limit theorem is provided here as well. These estimators and central limit theorems are results of this papers. They are applied to each of the 3 examples that are proposed in Section 2. In section 4 we provide a simulation study of parameter estimators for the considered copula families. For each of the proposed copulas, we prove existence of the MLE, and derive its asymptotic distribution. A comparative study is provided at the end of this section to show advantages of our proposed estimators. In Section 5 we provide comments and conclusions of the work. An appendix of proofs concludes the paper.

\section{Copula families of interest}
Copulas in general are defined as multivariate distributions with uniform marginals on (0,1) (see Nelsen (2006), Durante and Sempi (2016) or Sklar (1959) for more on the topic). Several authors have used this tool over the years to model dependence in various applied fields of research. In this paper, we are interested in copulas that were proposed by Longla (2023). Due to their novelty, these copulas have not yet been considered for estimation problems. There is no work that suggests how to use them, how to estimate their parameters or test any related claims.
Longla (2023) has shown that for any sequences $\lambda_k$ and $\varphi_k(x)$ such that: 
\begin{itemize}
\item The sequence of $\lvert{}\lambda_k\lvert{}$ has a finite number of values or converges to $0$.
\item The set $\{\varphi_k(x), k\in\mathbb{N}\}$ is an orthonormal basis  of $L^2(0,1)$.
\end{itemize} 
The function 
\begin{equation} c(u,v)=\sum_{k=1}^{\infty}\lambda_k\varphi_k(u)\varphi_k(v) \label{cop}
\end{equation}
is a copula density and the series converges pointwise and uniformly on $(0,1)^2$, if \begin{equation}
1+\sum_{k=1}^{\infty} \lambda_k \alpha_k \ge 0, \quad \text{where}\quad \alpha_k=\begin{cases} 
       \max\varphi^2_k & \text{if } \quad \lambda_k < 0 \\
      \min\varphi_k\max\varphi_k &\text{if} \quad \lambda_k > 0 .
   \end{cases} \label{cond2}
\end{equation}
Among examples of copula families of the form \eqref{cop}, that satisfy condition \eqref{cond2}, Longla (2023) introduced the following copulas.

\begin{itemize}
\item \textbf{The sine-cosine copulas as perturbations of $\Pi(u,v)$}
\end{itemize}
Longla (2023) considered the trigonometric orthonormal basis of $L^2(0,1)$ that consists of functions: $\{1, \sqrt{2}\cos 2\pi kx, \sqrt{2}\sin 2\pi kx, k\in\mathbb{N^*}\}$. He showed that 
\begin{eqnarray} \label{copt}
C(u,v)=uv+\frac{1}{2\pi^2}\sum_{k=1}^\infty \frac{\lambda_k}{k^2} \sin 2\pi ku \times \sin 2\pi kv+\nonumber\\ 
+\frac{1}{2\pi^2}\sum_{k=1}^\infty \frac{\mu_k}{k^2}(1-\cos 2\pi ku)(1-\cos 2\pi kv)
\end{eqnarray}
is a copula, when $\displaystyle \sum_{k=1}^{\infty} (\lvert{}\lambda_k\lvert{}+\lvert{}\mu_k\lvert{})\le 1/2$. This copula was called sine-cosine copula. It is a perturbation of the independence copula. 
 \begin{figure}[h!]
\begin{minipage}[c]{.46\linewidth}
     \begin{center}
             \includegraphics[width=4cm, height=4cm]{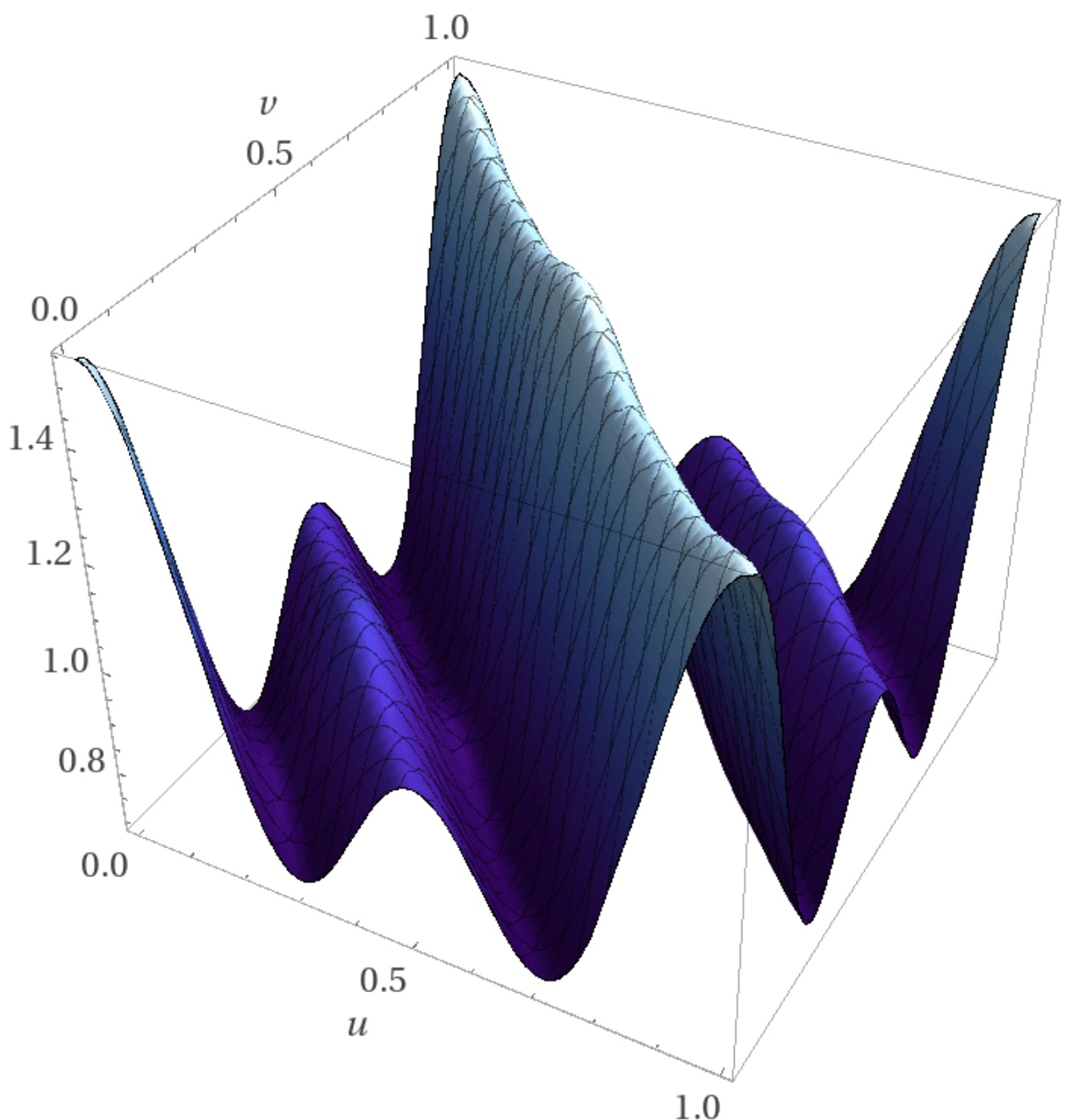}
         \end{center}
   \end{minipage} \hfill
   \begin{minipage}[c]{.46\linewidth}
    \begin{center}
            \includegraphics[width=4cm, height=4cm]{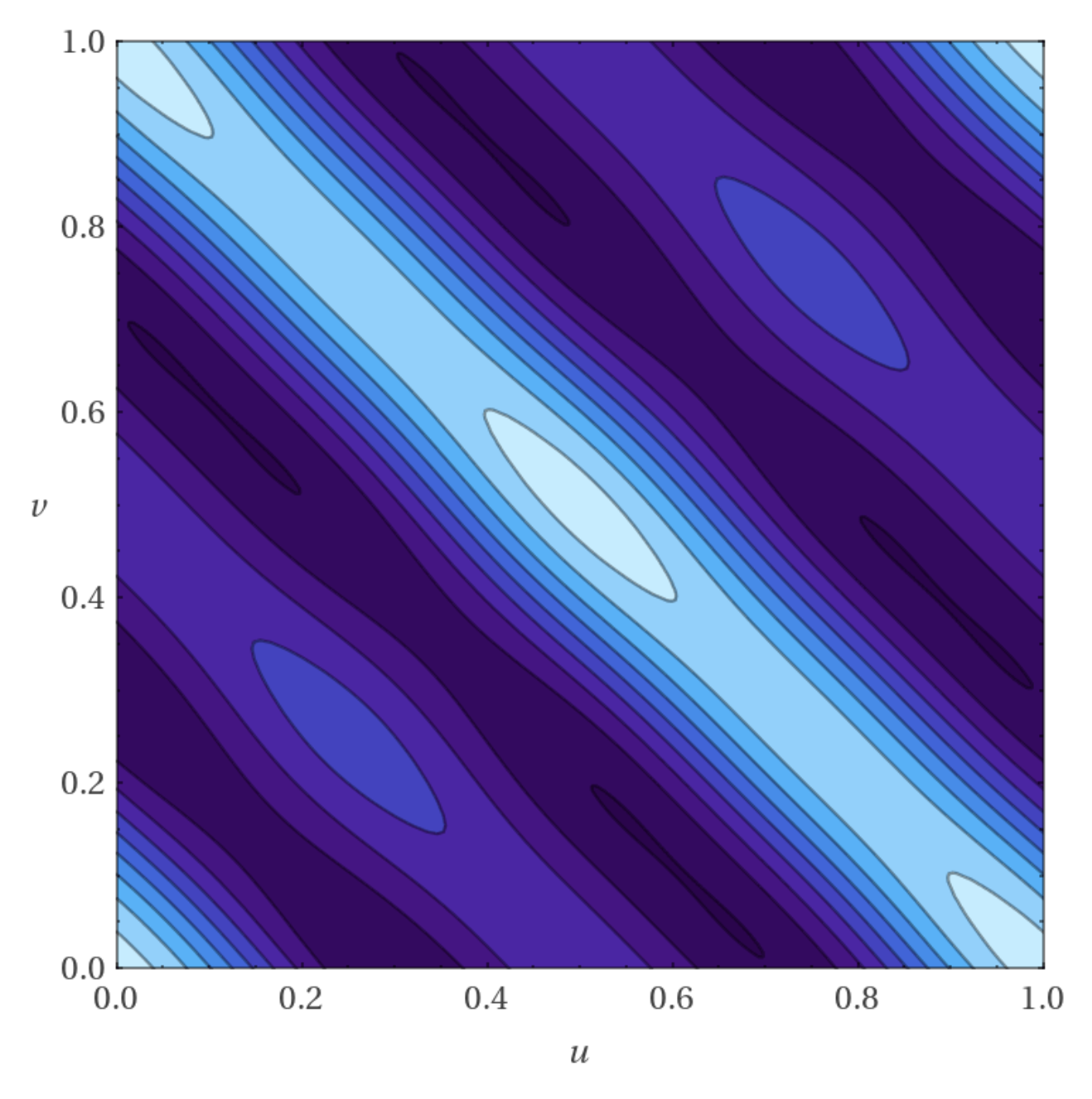}
        \end{center}
 \end{minipage}
 \caption{Sine-cosine copula: $\lambda_1=0.14$, $\mu_1=-0.11$ $\lambda_2=0.13$ and $\mu_2=-0.12$}
 \end{figure}
\begin{itemize}
\item \textbf{The sine copulas as perturbations of $\Pi(u,v)$}
\end{itemize}
Sine copulas were obtained using the basis of $L^2(0,1)$ given by $\{1, \sqrt{2}\cos k\pi x, k>0,k\in\mathbb{N} \}$. Their formula is
\begin{figure}[H]
\begin{minipage}[c]{.46\linewidth}
     \begin{center}
             \includegraphics[width=4cm, height=4cm]{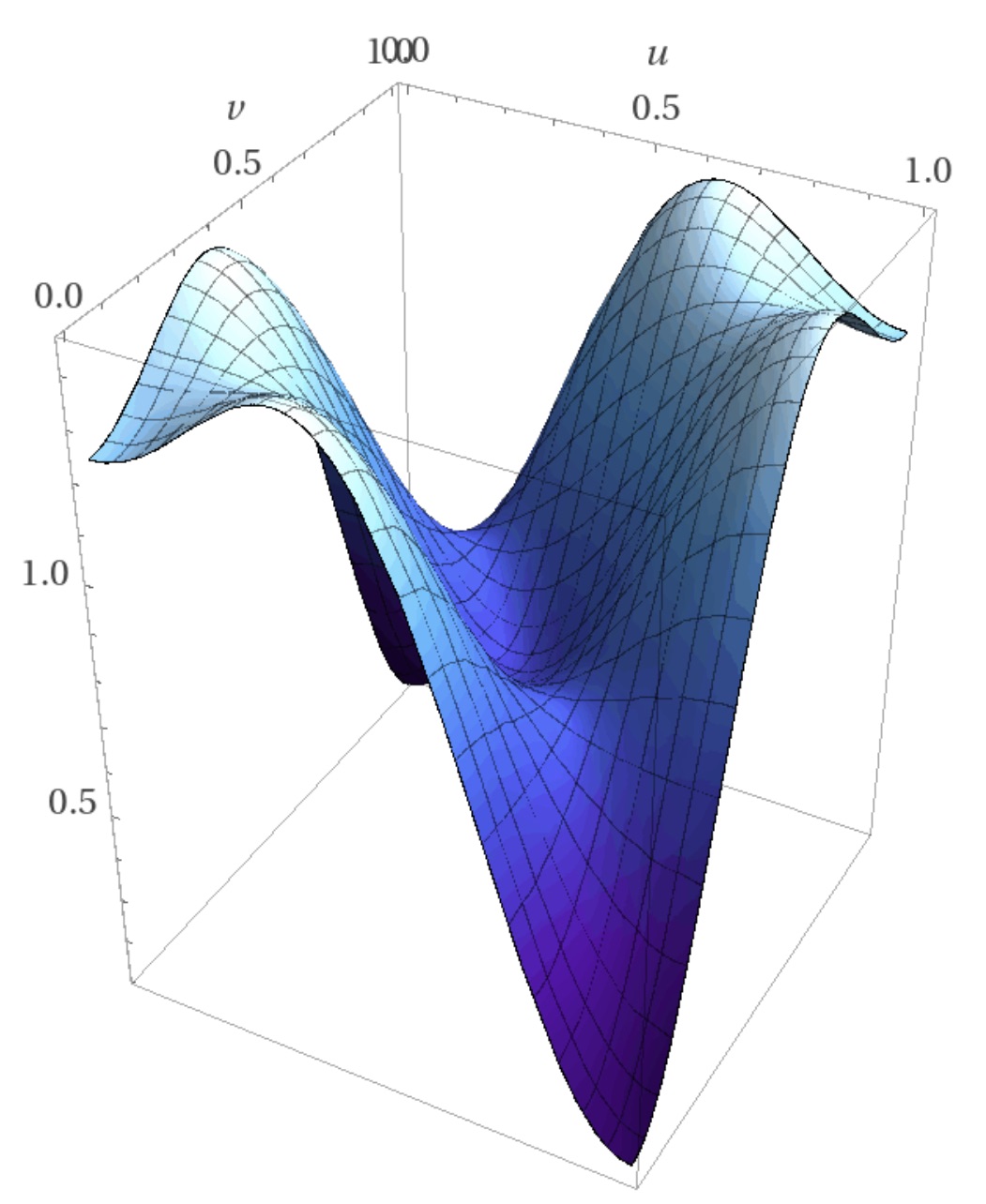}
         \end{center}
   \end{minipage} \hfill
   \begin{minipage}[c]{.46\linewidth}
    \begin{center}
            \includegraphics[width=4cm, height=4cm]{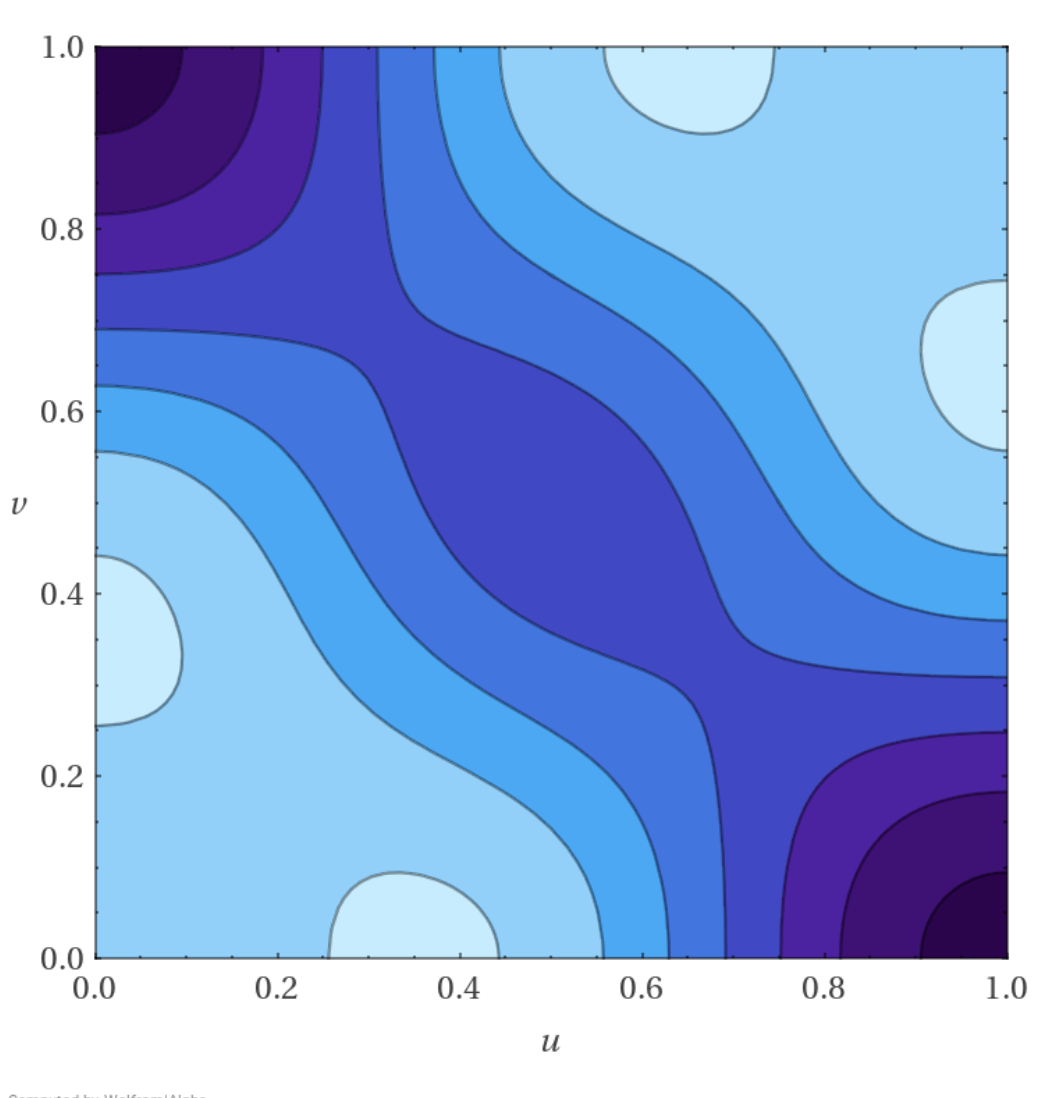}
        \end{center}
 \end{minipage}
 \caption{Sine copula with $\lambda_1=0.28$ and $\lambda_2=-0.15$}
 \end{figure}
\begin{equation}
C(u,v)=uv+\frac{2}{\pi^2}\sum_{k=1}^{\infty} \frac{\lambda_k}{k^2} \sin k\pi u \sin k\pi v, \quad \text{for} \label{cosG}
\quad
 \sum_{k=1}^\infty \lvert{}\lambda_k\lvert{}\leq 1/2.
\end{equation}
Any finite sum of terms from \eqref{cosG} that includes the first term $uv$  is a copula. Sine copulas and Sine-cosine copulas are not part of the trigonometric copulas of Chesneau (2021). These copulas are based on fundamental properties of linear operators induced by copulas on the Hilbert space $L^2(0,1)$. Their coefficients are eigen-values and the functions $\varphi_k(x)$ are eigen-functions of the operator (see Longla (2023)). 
 
\begin{itemize}
\item \textbf{The extended Farlie - Gumbel - Morgenstern copulas or Legendre copulas.}
\end{itemize}

They were obtained using shifted Legendre polynomials, orthonormal basis of $L^2(0,1)$ for $\{ 1, \varphi_k, k>0, k\in\mathbb{N}\}$. Legendre polynomials are defined by Rodrigues' formula  (see Rodrigues (1816))
\begin{equation}
P_k(x)=\frac{1}{2^k k!}\frac{d^k}{dx^k}(x^2-1)^k.
\end{equation}

 \begin{figure}[H]
\begin{minipage}[c]{.46\linewidth}
     \begin{center}
             \includegraphics[width=4cm, height=4cm]{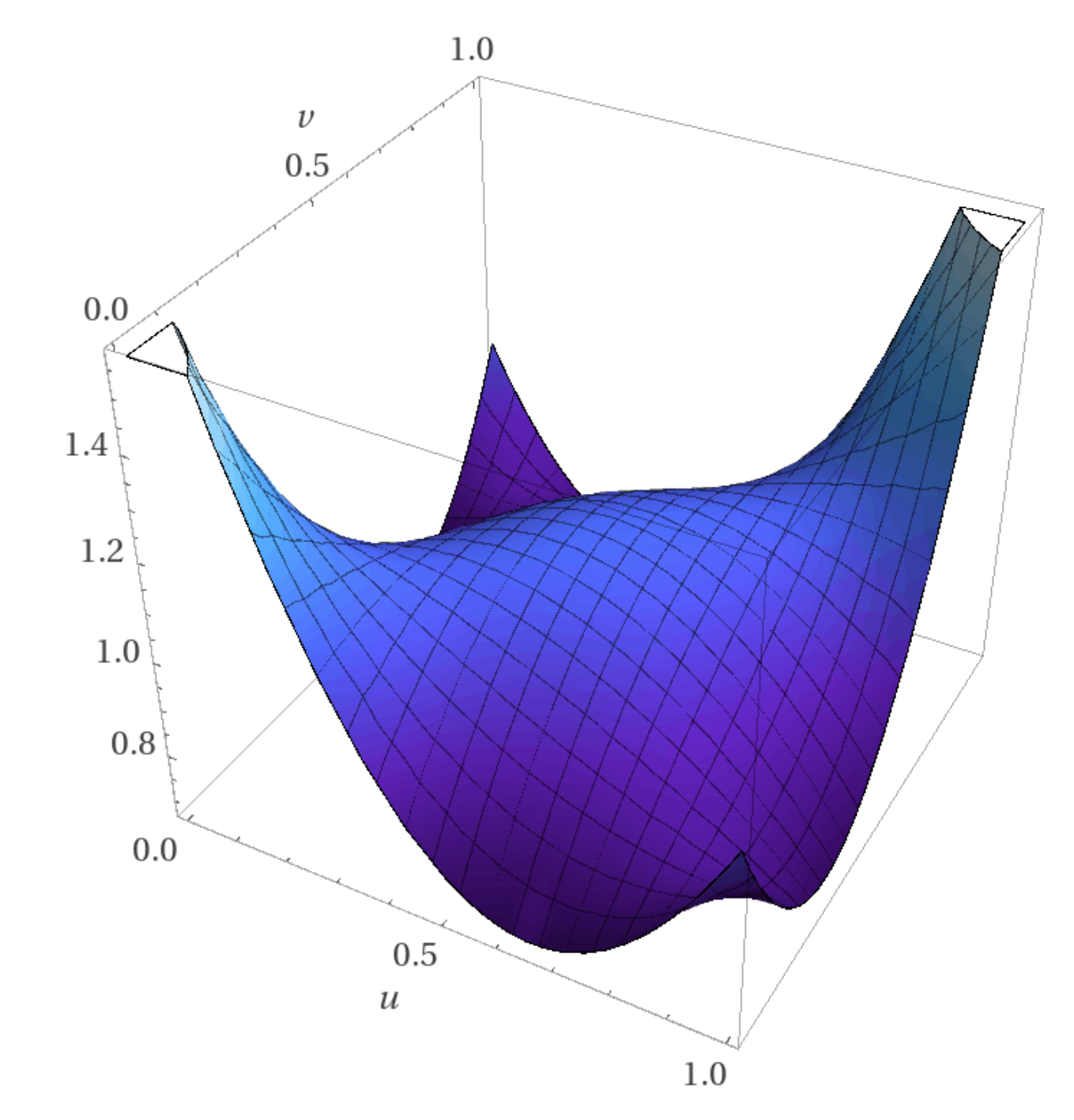}
         \end{center}
   \end{minipage} \hfill
   \begin{minipage}[c]{.46\linewidth}
    \begin{center}
            \includegraphics[width=4cm, height=4cm]{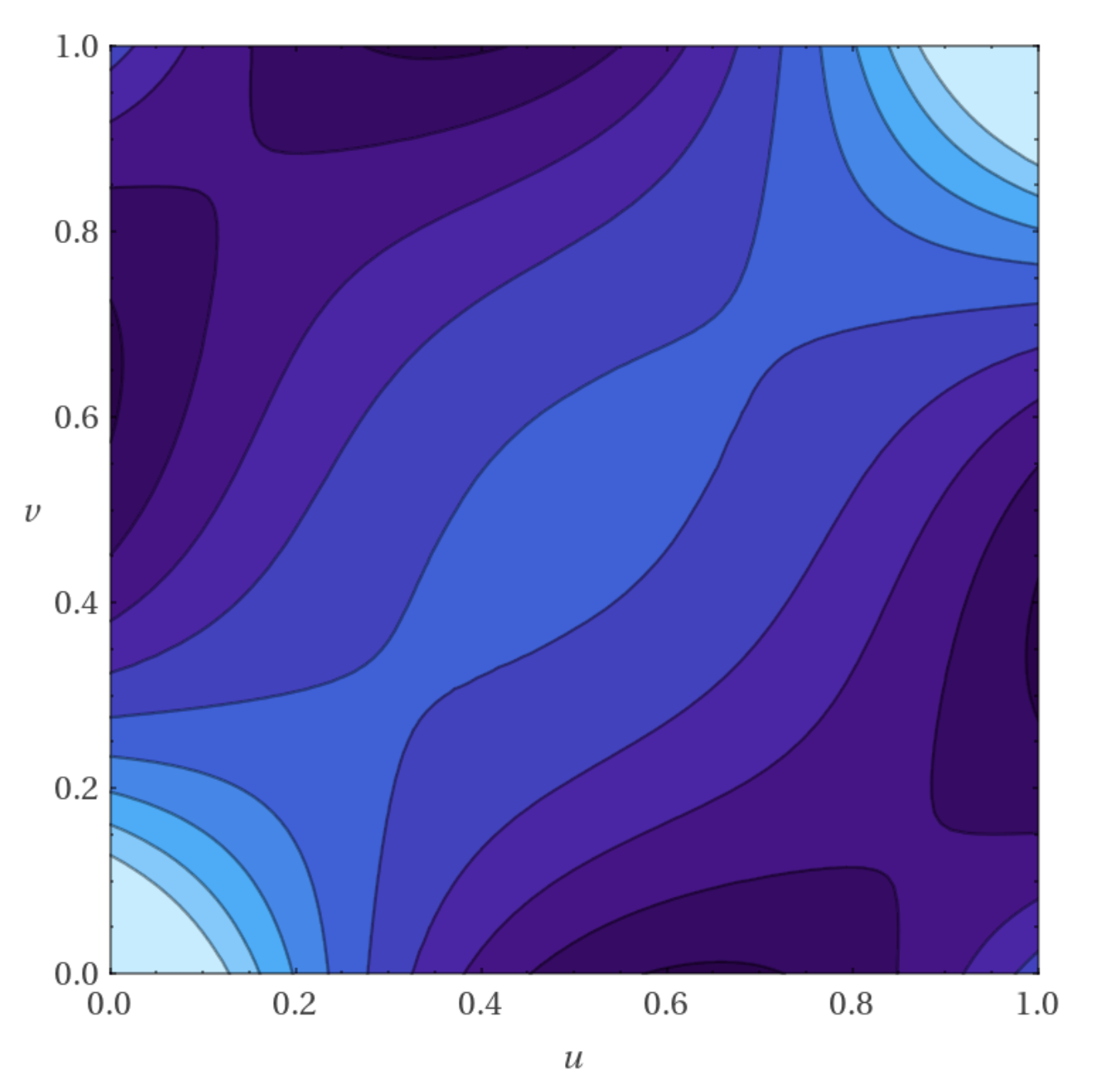}
        \end{center}
 \end{minipage}
 \caption{Legendre copula with $\lambda_1=0.15$ and $\lambda_2=0.1$}
 \end{figure}

The orthonormal basis of shifted Legendre polynomials is
\begin{equation}
\varphi_k(x)=\sqrt{2k+1}P_k(2x-1).
\end{equation}

For $k>0,$ $\varphi_k(1)=\sqrt{2k+1}=\max\varphi_k$, for odd  $k$, $\min \varphi_k=-\sqrt{2k+1}=\varphi_k(0)$ and for even $k,$ $0>\min\varphi_k>-\sqrt{2k+1}$. This minimum is achieved at least at one point. Therefore, when \eqref{cop} contains $\varphi_k(x)$ alone, the range of $\lambda_k$ is \[ \frac{-1}{2k+1}\le \lambda_k\leq \min(1, \frac{1}{(2k+1)\lvert{}\min P_k(x)\lvert{}}). \] 

Using this basis in formula \eqref{cop}, we obtain what Longla (2023) called extended FGM copula family for any sum of terms that contains the first two functions. For more on Legendre polynomials, see Belousov (1962), Szeg\"{o} (1975) and the references therein.

 \section{Parameter estimation }
The constructed copula families include the extended Farlie-Gumbel-Morgenstern copula family. The FGM copula family has been extensively studied in recent years. Earlier studies include Morgenstern (1956), Farlie (1960) and Gumbel (1960). These authors discussed  the bivariate FGM copula family. Johnson and Kotz (1975) formulated the FGM copula as a multivariate distribution on $[0,1]^d$ for any positive integer $d>1$. The multivariate FGM copula family is just as valuable as  the multivariate normal distribution thanks to its simplicity. It has been applied to statistical modeling in various research fields such as economics (see Patton (2006)), finance (see Cossette et al. (2013)) and reliability engineering (see Navarro and Durante (2017) and recently Ota and Kimura (2021)).

Suppose that we have a Markov chain generated by a copula and the uniform distribution as stationary distribution.
Assuming that this copula has the representation \eqref{cop}, we have the following for estimators of the copula coefficients.
\begin{theorem} \label{lambda1}
For any copula-based Markov chain $(U_0,U_1, \cdots, U_n)$ generated by the copula \eqref{cop} and the uniform marginal distribution,
\begin{equation}
\hat{\lambda}_k=\frac{1}{n}\sum_{i=1}^{n}\varphi_k(U_i)\varphi_k(U_{i-1})
\end{equation}
is an unbiased consistent estimator of $\hat{\lambda}_k$ provided that $$\sum \lvert{}\lambda_z\lvert{} (\int_0^1\varphi_z(x)\varphi_k^2(x)dx)^2)< \infty \quad \text{and} \quad  \int_0^1 \varphi_k^4(x)dx<\infty.$$
Moreover,    $\sqrt{n}(\hat{\lambda}_k-\lambda_k)\to N(0,\sigma^2)$, where the positve real number $\sigma^2$ is given by  
$$\sigma^2=(1-\lambda_k^2+\sum \lambda_z (\int_0^1\varphi_z(x)\varphi_k^2(x)dx)^2)+$$ $$+2\lambda_k^2(\int_0^1\varphi^4_k(x)dx-1)+2\sum_{m=2}^{\infty}\lambda_k^2\sum\lambda^{m-1}_z(\int_0^1\varphi_z(x)\varphi_k^2(x)dx)^2.$$
\end{theorem}

Note, that when $\varphi_k(x)$ are uniformly bounded, the conditions on $\lambda_k$ reduce to $\sum \lvert \lambda_k \lvert <\infty,$ which is satisfied under condition \eqref{cond2}. 
Based Theorem \eqref{lambda1}, we can construct confidence intervals for each of the values of $\lambda_k$ separately for every case of $\varphi_k(x)$. This implies estimation of the variance as shown in the examples below. But, to tackle any problem related to joint parameters, we need to establish a result for the joint limiting distribution of these estimators of the coefficients of our copula. These distributions can be valuable in model selection for the number of parameters to include in the estimator of the copula. We need to find the distribution of $\Lambda=(\hat{\lambda}_{k_1}, \cdots, \hat{\lambda}_{k_s})$ for any sequence of increasing integers $\hat{\lambda}_{k_1}, \cdots, \hat{\lambda}_{k_s}$.

\begin{theorem} \label{lambdam}
Assume $(U_0,U_1, \cdots, U_n)$ is a copula-based Markov chain generated by a copula of the form \eqref{cop}. Assume $\int_{0}^1\varphi^2_{k_i}(x)\varphi^2_{k_j}(x)dx<\infty$, $$\sum_{z=1}^{\infty}\lvert{}\lambda_z\lvert{} (\int_0^1\varphi_{z}(x)\varphi_{k_i}(x)\varphi_{k_j}(x)dx)^2<\infty \quad \text{and}$$
$$\sum_{z=1}^{\infty}\lvert{}\lambda_z\lvert{} (\int_0^1\varphi_{z}(x)\varphi_{k_i}^2(x)dx)(\int_0^1\varphi_{z}(x)\varphi_{k_j}^2(x)dx)<\infty,$$ then for some positive definite matrix  $\Sigma$ we have  $$\sqrt{n}(\hat{\Lambda}-\Lambda) \to N_s({\bf 0}, \Sigma).$$ 
\end{theorem}

\begin{remark}
For an orthonormal basis of $L^2(0,1)$ made of uniformly bounded and symmetric functions, the conditions of Theorem \ref{lambdam} are automatically satisfied because the integrals are uniformly bounded and $\sum_{i=1}^\infty \lvert{}\lambda_i\lvert{}<\infty$. 
\end{remark}

\begin{proposition}\label{prop1}
    Under the assumptions of Theorem \ref{lambdam}, the off-diagonal elements of the matrix $\Sigma$ are given by
    \begin{eqnarray}\label{cov}
\underset{k\neq j}{\Sigma_{k,j}}&=&\sum\limits_{z=1}^s\lambda_z(\int_0^1\varphi_z(x)\varphi_k(x)\varphi_j(x)dx)^2+\lambda_k\lambda_j\left(2(\int_0^1\varphi^2_k(x)\varphi_j^2(x)dx-1)\right.\nonumber\\&&\left.+2\sum\limits_{z=1}^s\sum\limits_{m=2}^\infty\lambda_z^{m-1}(\int_0^1\varphi_z(x)\varphi_k^2(x)dx)(\int_0^1\varphi_z(x)\varphi_j^2(x)dx)-1\right).
\end{eqnarray}
\end{proposition}

\subsection{Applications to some examples}
For the classes of copulas given above, the Spearman's $\rho(C)$ is estimated by  $\hat{\lambda}_1$.
\begin{enumerate}
\item Sine-cosine copulas. Theorem \ref{lambda1} implies that for the copulas based on the sine-cosine basis, 
$$\displaystyle \sigma^2_k=1+\frac{\lambda_{2k}}{2}+\frac{\lambda_k^2\lambda_{2k}}{1-\lambda_{2k}}, \quad \text{for} \quad  \hat{\lambda}_k=\frac{1}{n}\sum_{i=1}^{n} 2\cos 2\pi k U_i \cos 2\pi k U_{i-1}$$
 $$\text{and} \quad  \sigma^2_k=1+\frac{\lambda_{2k}}{2}+\frac{\mu_k^2\lambda_{2k}}{1-\lambda_{2k}}, \quad \text{for}\quad  \hat{\mu}_k=\frac{1}{n}\sum_{i=1}^{n} 2\sin 2\pi k U_i \sin 2\pi k U_{i-1}.$$

It is also clear that for sine-cosine copulas, with $i\ne j$, $\Sigma_{i,j}$ is given by

\begin{center}
\begin{tabular}{ | c | c| c | } 
  \hline
$\varphi_{k_i}/\varphi_{k_j}$ & $\varphi_{1}(x)=\sqrt{2}\cos2\pi k_j x$ & $\varphi_{2}(x)=\sqrt{2}\sin 2\pi k_j x$ \\ 
  \hline
 $\varphi_1(x)$ & $\frac{1}{2}(\lambda_{k_i+k_j}+\lambda_{\lvert{}k_i-k_j\lvert{}})-\lambda_{k_i}\lambda_{k_j}$ & $\frac{1}{2}(\mu_{k_i+k_j}+\mu_{\lvert{}k_i-k_j\lvert{}})-\lambda_{k_i}\mu_{k_j}$ \\ 
  \hline
  $\varphi_2(x)$ & $\frac{1}{2}(\mu_{k_i+k_j}+\mu_{\lvert{}k_i-k_j\lvert{}})-\mu_{k_i}\lambda_{k_j}$ & $\frac{1}{2}(\lambda_{k_i+k_j}+\lambda_{\lvert{}k_i-k_j\lvert{}})-\mu_{k_i}\mu_{k_j}$ \\ 
  \hline
\end{tabular}
\end{center}

Therefore, the multivariate central limit theorem holds with $\Sigma_{ij}$ as entries of the covariance matrix $\Sigma$.
A test of independence can be formed based on this limit theorem. Under the hypothesis of independence of observations agains that of a  Markov chain with copula $C(u,v)$, all coefficients are equal to $0$. The estimators are asymptotically independent with $\sigma^2=1$. We can formulate this in the following proposition.
\begin{proposition}\label{lambdacosine}
For any square integrable copula with representation in sine-cosine series, under the assumption that all parameters are equal to zero, any vector of estimators with finitely many components satisfies the $\sqrt{n}$-central limit theorem with variance-covariance matrix $\Sigma=\mathbb{I}_s$, where $s$ is the length of the vector and $\mathbb{I}_s$ is the identity matrix of dimension $s$. Moreover, if a sequence $W_{n_i}$ takes on values of $\hat{\lambda}_{k}$ or $\hat{\mu}_{k}$ for the Markov chain $U_0,\cdots,U_n$, then 
$$n\sum_{i=1}^s W_{n_i}^2\to \chi^2(s)\quad  \text{as} \quad n\to\infty \quad \text{and}\quad \frac{1}{\sqrt{2s}}(n\sum_{i=1}^s W_{n_i}^2-s)\to N(0,1)$$ $$ \text{as}\quad  n\to\infty,\quad \text{then}\quad s\to\infty.$$
\end{proposition}

\item Sine copulas.
Considering  sine copulas, we have $\displaystyle \sigma^2_k=1+\frac{\lambda_{2k}}{2}+\frac{\lambda_k^2\lambda_{2k}}{1-\lambda_{2k}}.$ Estimators are $$ \hat{\lambda}_k=\frac{1}{n}\sum_{i=1}^{n} 2\cos \pi k U_i \cos \pi k U_{i-1}.$$ Therefore, Proposition \ref{lambdacosine} holds exactly with the same conditions.
\end{enumerate}

\section{Simulation study for the 3 examples}
In this section, we compare these estimators to other known estimators such as maximum likelihood estimator, the robust estimator of Longla and Peligrad (2021) in two forms. The performance is studied on Markov chains simulated using three different copulas from the families described above. 

\subsection{ The specific copula examples}
We consider 3 types of copulas with densities of the form \eqref{cop}. A sine copula, a sine-cosine copula and a Legendre copula. The first two present several tops and bottoms, while the second is an extension of the FGM copula family. This third copula can be used to test departure from the FGM copula family in modelling. They are defined as follows.
\begin{enumerate}
\item {\it Sine copulas}\\
 \begin{equation}\label{e2}
C(u,v)=uv+\dfrac{2}{\pi^2}\lambda_1 \sin(\pi u)\sin(\pi v)+\dfrac{1}{2\pi^2}\lambda_2\sin(2\pi u)\sin(2\pi v),\end{equation}
for $\lvert \lambda_1\lvert+\lvert \lambda_2\lvert \le 0.5$ and density 
\begin{equation*} \label{e3}
c(u,v)=1+2\lambda_1 \cos(\pi u)\cos(\pi v)+2\lambda_2 \cos(2\pi u)\cos(2\pi v).\end{equation*}
\item {\it Sine-cosine copulas}\\
\begin{eqnarray*}\label{e4}
C(u,v)&=&uv+\dfrac{1}{2\pi^2}(\lambda_1\sin(2\pi u)\sin(2\pi v)+\mu_1(1-\cos(2\pi u))(1-\cos(2\pi v))\nonumber\\&&+\dfrac{1}{8\pi^2}(\lambda_2\sin(4\pi u)\sin(4\pi v)+\mu_2(1-\cos(4\pi u))(1-\cos(4\pi v)), \end{eqnarray*}
for $\lvert \lambda_1\lvert+\lvert \lambda_2\lvert+\lvert \mu_1\lvert+\lvert \mu_2\lvert\le 0.5$ and density \begin{eqnarray}\label{e5}
c(u,v)&=&1+2\lambda_1\cos(2\pi u)\cos(2\pi v)+2\mu_1\sin(2\pi u)\sin(2\pi v)\nonumber\\&&+2\lambda_2\cos(4\pi u)\cos(4\pi v)+2\mu_2\sin(4\pi u)\sin(4\pi v)).\end{eqnarray}
\item {\it Legendre copulas}\\
If $\varphi_k\in\{\sqrt{3}(2x-1),~ \sqrt{5}(6x^2-6x+1)\}$ and $ 3\vert \lambda_1\vert+5\vert \lambda_2\vert\leq 1$, we obtain the following Legendre copula: \begin{equation}\label{e8} C(u,v)=1+3\lambda_1(u^2-v)(v^2-v)+5\lambda_2(2u^3-3u^2+u)(2v^3-3v^2+v),\end{equation} with density
\begin{equation}\label{e6}
c(u,v)=1+3\lambda_1(2u-1)(2v-1)+5\lambda_2(6u^2-6u+1)(6v^2-6v+1).\end{equation} \end{enumerate}
Each of these copulas is a perturbation of the independence copula and can be used to seek some level of dependence in the data. We generate Markov chains using copulas from these families for known values of parameters. Then, we use them to build confidence intervals, and compute coverage probabilities based on the asymptotic distrubtions of estimators that are proposed.

\subsection{Algorithm for simulation of the Markov chains}\label{sec2}
As shown in Darsow et al. (1992), the transition kernel for a copula-based Markov chain is given by the derivative with respect to the first variable of the copula by $P(x, [0,y])=C_{,1}(x,y)$. Therefore, to generate a realization of a Markov chain $(U_0,\cdots, U_n)$ using each of the three considered copulas and the uniform marginal distribution, we use the following algorithm.
 \begin{enumerate}
\item generate $U_0 \sim Unif(0,1)$;
\item $U[1]=U_0$;
\item For $2\leq i\leq n+1$: 
 
a). generate $w\sim Unif(0,1)$;
 
b). $U[i]$ is the unique root on $[0,1]$ of the equation $ C_{,1}(U[i-1],U[i])-w=0$. 
\end{enumerate}

For each of the used copulas, the last step uses the $R$ function $Uniroot$.  For n=999, we have obtained the following graphs, displaying time series {\it copula 1, copula 2} and {\it copula 3} respectively for {\it sine copula, sine-cosine copula} and {\it Legendre copula}.

\begin{figure}[H]\begin{center}\includegraphics[width=13cm, height=7cm]{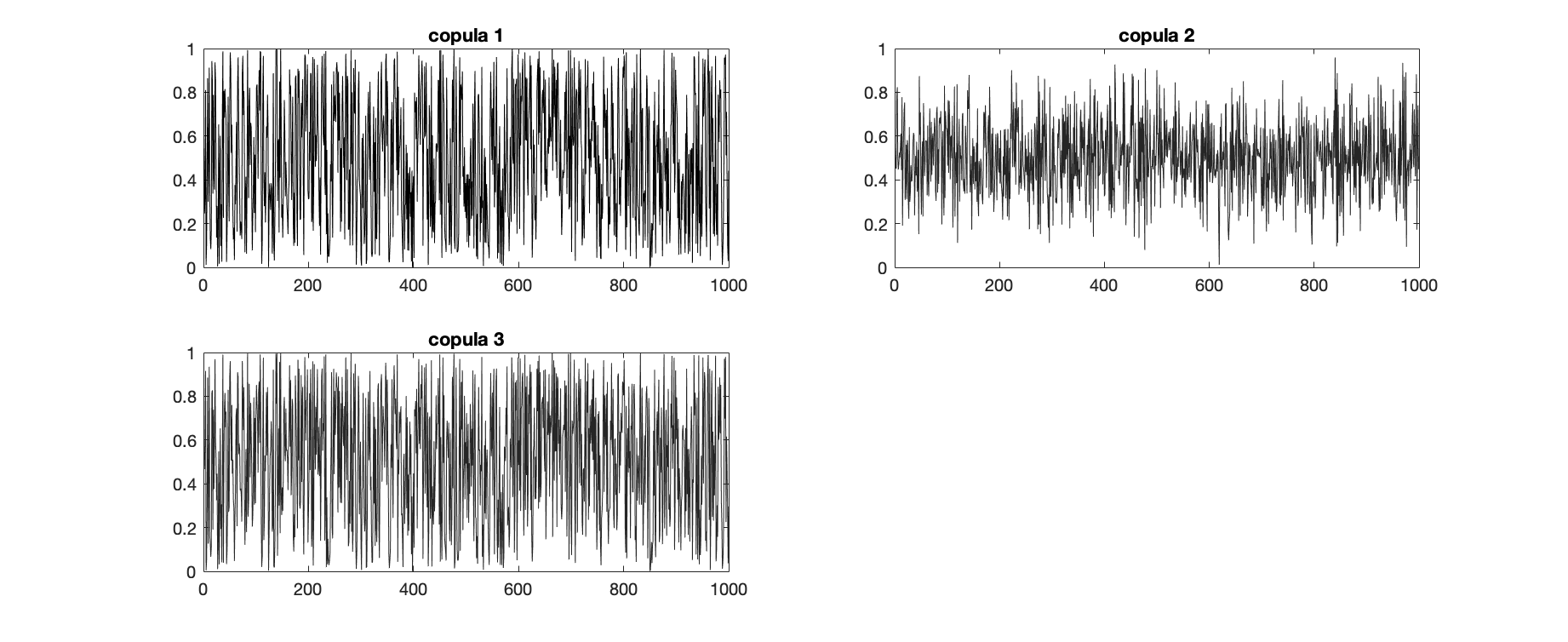}\caption{\it Data generated under  each of the copula for n=999. }\label{ill}\end{center}\end{figure} 

\subsection{Study of the proposed estimator for our examples}
Assume a Markov chain $(U_0,\cdots, U_n)$ is generated using a copula of the form \eqref{cop}.
The estimator of $\lambda_k$ is given by
\begin{equation*}\label{1}\hat{\lambda}_k=\dfrac{1}{n}\sum\limits_{i=1}^{n}\varphi_k(U_i)\varphi_k(U_{i-1}),\end{equation*}
and $\hat{\Lambda}=(\hat{\lambda}_1,...,\hat{\lambda}_s)$,  satisfies
$
\sqrt{n}(\hat{\Lambda}-\Lambda)\longrightarrow {N}\left({\bf 0}, \Sigma\right),$ where $\Sigma$ is the variance-covariance matrix defined in the Proposition \ref{components} below. 

\begin{proposition}\label{components} Copulas of the examples are constructed by means of uniformly bounded functions $\varphi_k(x)$ and the sum of absolute values of parameters is bounded. Therefore, Theorem \ref{lambdam} holds with the following $\Sigma$. 
\begin{enumerate}
\item \noindent For the sine copula of formula (\ref{e2}) with $\hat{\Lambda}=(\hat{\lambda}_1,\hat{\lambda}_2)$:\begin{equation}\label{matrix1}
    \Sigma=\left(\begin{matrix}
    1+\dfrac{\lambda_{2}}{2}+\dfrac{\lambda_1^2\lambda_{2}}{1-\lambda_{2}}&~~ \lambda_1\left(\dfrac{1}{2}-\lambda_2\right)\\~\\
    \lambda_1\left(\dfrac{1}{2}-\lambda_2\right)&1
    \end{matrix}\right).
\end{equation}
\item For the Legendre copula of formula (\ref{e8}) with $\hat{\Lambda}=(\hat{\lambda}_1,\hat{\lambda}_2)$:
\begin{equation}\label{matrix3}
    \Sigma=\left(\begin{matrix}
    1+\dfrac{4}{5}\lambda_2+\lambda_1^2\dfrac{3+5\lambda_2}{5(1-\lambda_2)}&~~ \dfrac{4}{5}\lambda_1+\lambda_1\lambda_2\dfrac{1+7\lambda_2}{7(1-\lambda_2)}\\~\\
   \dfrac{4}{5}\lambda_1+\lambda_1\lambda_2\dfrac{1+7\lambda_2}{7(1-\lambda_2)} &~~1+\dfrac{20}{49}\lambda_2+\lambda_2^2\dfrac{63-23\lambda_2}{49(1-\lambda_2)}
    \end{matrix}\right).
\end{equation}

\item \noindent For the sine-cosine copula with density (\ref{e5}) and $\hat{\Lambda}=(\hat{\lambda}_1,\hat{\lambda}_2,\hat{\mu}_1,\hat{\mu}_2)$:
 \begin{equation}\label{matrix2}
   \left(\begin{matrix}
    1+\dfrac{\lambda_{2}}{2}+\dfrac{\lambda_1^2\lambda_{2}}{1-\lambda_{2}}&~~ \lambda_1\left(\dfrac{1}{2}-\lambda_2\right)&~~\dfrac{\mu_2}{2}-2\lambda_1\mu_1& ~~ \dfrac{\mu_1}{2}-\lambda_1\mu_2\\~\\
    \lambda_1\left(\dfrac{1}{2}-\lambda_2\right)&~~ 1&~~\dfrac{\mu_1}{2}-\lambda_2\mu_1&~~ -2\lambda_2\mu_2\\~\\
    \dfrac{\mu_2}{2}-2\lambda_1\mu_1&~~\dfrac{\mu_1}{2}-\lambda_2\mu_1&~~ 1+\dfrac{\lambda_{2}}{2}+\dfrac{\mu_1^2\lambda_{2}}{1-\lambda_{2}}&~~\dfrac{\lambda_1}{2}-\mu_1\mu_2\\~\\
   \dfrac{\mu_1}{2}-\lambda_1\mu_2&~~-2\lambda_2\mu_2&~~\dfrac{\lambda_1}{2}-\mu_1\mu_2 &~~ 1
    \end{matrix}\right).
\end{equation}

\end{enumerate}
\end{proposition}

Proposition \ref{components} shows that the estimators are asymptotically independent if $\lambda_1=0$ or $\lambda_2=0.5$ for the sine copula. This is true for the Legendre copula if $\lambda_1=0$. These matrices allow us to construct confidence intervals for single copula parameters. Any test of hypotheses on these parameters can use the obtained multivariate central limit theorem.

\begin{remark}\label{rmk1}
Under the assumptions of Proposition \ref{components}, the individual confidence interval for each of the parameters is given by \begin{equation}\label{ci1} \left(\hat{\lambda}_k-z_{\frac{\alpha}{2}}\frac{\hat{\sigma}_k}{\sqrt{n}},\quad  \hat{\lambda}_k+z_{\frac{\alpha}{2}}\frac{\hat{\sigma}_k}{\sqrt{n}}\right)\text{ and }\left(\hat{\mu}_k-z_{\frac{\alpha}{2}}\frac{\hat{\tau}_k}{\sqrt{n}},\quad \hat{\mu}_k+z_{\frac{\alpha}{2}}\frac{\hat{\tau}_k}{\sqrt{n}}\right),\end{equation} 
\[ \text{Moreover,}\quad  Q=n(\hat{\Lambda}-\Lambda)'\Sigma^{-1}(\hat{\Lambda}-\Lambda)\underset{n\rightarrow \infty}{\longrightarrow} \chi^2_s, \quad \text{and}\] 
\begin{eqnarray}C=\left\{\Lambda\in\mathbb{R}^s: \hat{Q}< \chi^2_\alpha(s) \right\}
=\left\{\Lambda\in\mathbb{R}^s: n(\hat{\Lambda}-\Lambda)'\hat{\Sigma}^{-1}(\hat{\Lambda}-\Lambda)< \chi_\alpha^2(s),\right\}
\end{eqnarray}

is the $(1-\alpha)100\%$ confidence region for the vector $\Lambda$,
$\chi^2_{\alpha}(s)$ is given by $P(\chi^2(s)>\chi^2_\alpha(s))=\alpha$,  
 $\hat{\sigma}_k$ ($\hat{\tau}_k$) is the estimate of $\sigma_k$ ($\tau_k$) obtained by replacing $\lambda_k$ and $\mu_k$ by their estimates in the matrix $\Sigma$. The last comment is true thanks to consistency of the estimators.
\end{remark}

\subsection{The Maximum likelihood estimator (MLE)}
For the Markov chain $(U_0, \cdots, U_n)$ generated by each of the considered copulas and the uniform marginal distribution, the log-likelihood function is:
\begin{enumerate}
\item For the sine copula:
\[\ell(\Lambda)=\sum\limits_{i=1}^n \ln(1+2\lambda_1 \cos(\pi U_i)\cos(\pi U_{i-1})+2\lambda_2 \cos(2\pi U_i)\cos(2\pi U_{i-1})),\]
 and  the MLE of $\Lambda=(\lambda_1,\lambda_2)$ is given by:
\begin{equation}
\hat{\Lambda}^{ML}=\underset{\lvert\lambda_1\lvert+\lvert\lambda_2\lvert\leq .5}{Argmax}\ell(\Lambda).
\end{equation}
\item For the sine-cosine copula:
\begin{eqnarray*}\ell(\Lambda)=\sum\limits_{i=1}^n \ln(1+2\lambda_1\cos(2\pi U_i)\cos(2\pi U_{i-1})+2\mu_1\sin(2\pi U_i)\sin(2\pi U_{i-1})\\+2\lambda_2\cos(4\pi U_i)\cos(4\pi U_{i-1})+2\mu_2\sin(4\pi U_i)\sin(4\pi U_{i-1})),\end{eqnarray*} and  the MLE of $\Lambda=(\lambda_1,\lambda_2,\mu_1,\mu_2)$  is given by:
\begin{equation}
\hat{\Lambda}^{ML}=\underset{\lvert\lambda_1\lvert+\lvert\lambda_2\lvert+\lvert\mu_1\lvert+\lvert\mu_2\lvert\leq .5}{Argmax}\ell(\Lambda).
\end{equation}
\item For the Legendre copula:
\[\ell=\sum\limits_{i=1}^n \ln(1+3\lambda_1(2U_i-1)(2U_{i-1}-1)+5\lambda_2(6U_i^2-6U_i+1)(6U_{i-1}^2-6U_{i-1}+1)) ,\] and  the MLE of $\Lambda$ is given by:
\begin{equation}
\hat{\Lambda}^{ML}=\underset{3\lvert\lambda_1\lvert+5\lvert\lambda_2\lvert\leq 1}{Argmax}\ell(\Lambda).
\end{equation}
\end{enumerate}

Billingsley (1961) proved asymptotic likelihood theory for the MLE of parameters of a Markov process ($\sqrt{n}(\hat{\Lambda}^{ML}-\Lambda)\to N({\bf 0}, \Sigma^{-1})$) under the following regularity conditions.
\begin{enumerate}
\item The Markov chain is ergodic and not necessarity stationary.
\item For $u\in[0,1]$,  the set of $v$ for which $c(u,v)>0$ does not depend on $\Lambda$.
\item For $1\leq i,j,k\leq s$, where $s$ is the length of the vector $\Lambda$,\\
    \[\frac{\partial}{\partial\lambda_i}\ln c(u,v),\frac{\partial^2}{\partial\lambda_i\lambda_j}\ln c(u,v)\text{ and }\frac{\partial^3}{\partial\lambda_i\lambda_j\lambda_k}\ln c(u,v) \] exist and are continuous on the set of all considered $\Lambda$.
\item For $1\leq i,j,k\leq s$, $A'$ an open subset of the support of the vector $\Lambda$:\\
$\mathbb{E}_\Lambda\left[\left(\frac{\partial}{\partial\lambda_i}\ln c(u,v)\right)^2\right]<\infty$ and $\mathbb{E}_\Lambda\left[\underset{\Lambda\in A'}{\sup}\lvert\frac{\partial^3}{\partial\lambda_i\lambda_j\lambda_k}\ln c(u,v)\lvert\right]<\infty$.
\item  \begin{equation}\label{matrixx}\sigma_{i,j}=\mathbb{E}_\Lambda\left[\frac{\partial}{\partial\lambda_i}\ln c(u,v)\frac{\partial}{\partial\lambda_j}\ln c(u,v)\right]\end{equation} form a non-singular matrix $\Sigma$.
\end{enumerate}

\begin{theorem}\label{MLE1}
For any Markov chain generated by the considered copulas and the uniform distribution on $[0,1]$, there exist a consistent MLE $\Lambda^{ML}$ such that
\[\sqrt{n}(\Lambda^{ML}-\Lambda)\to N({\bf 0}, \Sigma^{-1})\]
and the $(1-\alpha)100\%$ confidence intervals for $\lambda_k$ (or $\mu_k $ for the sine-cosine copula) is
\begin{equation}\label{ci2}
\left({\lambda}^{ML}_k-z_{\dfrac{\alpha}{2}}\sqrt{\left[I_n^{-1} ({\Lambda}^{ML})\right]_{k,k} }~,~  {\lambda}^{ML}_k+z_{\dfrac{\alpha}{2}}\sqrt{ [I_n^{-1} ({\Lambda}^{ML})]_{k,k}} \right).
\end{equation}
\end{theorem}

\begin{remark} In each of the cases, we have the following.

1.  For the sine copula:\[ I_n(\Lambda)=\left(\sum_{i=1}^n\frac{4\cos(k\pi U_i)\cos(k\pi U_{i-1})\cos(j\pi U_i)\cos(j\pi U_{i-1})}{c^2(U_i,U_{i-1})}\right)_{1\leq k,j\leq 2}\] where $c(u,v)$ is the density of the cosine copula. 

2.  For the sine-cosine copula:
$\displaystyle
I_n(\Lambda)=\left(a_{j,k}, 
    \right)_{1\leq j,k\leq 4}\quad \text{where}$
~~~~~~ \begin{eqnarray*}
a_{1,1}=\sum_{i=1}^n\dfrac{4C^2_1(U_i)C^2_1(U_{i-1})}{c^2(U_i,U_{i-1})}&;& a_{1,2}=\sum_{i=1}^n\dfrac{4C_1(U_i)S_1(U_i)C_1(U_{i-1})S_1(U_{i-1})}{c^2(U_i,U_{i-1})};\end{eqnarray*} \begin{eqnarray*}
a_{1,3}=\sum\limits_{i=1}^n\dfrac{4C_1(U_i)C_2(U_{i})C_1(U_{i-1})C_2(U_{i-1})}{c^2(U_i,U_{i-1})} \end{eqnarray*} \begin{eqnarray*} a_{1,4}=\sum\limits_{i=1}^n\dfrac{4C_1(U_i)S_2(U_{i})C_1(U_{i-1})S_2(U_{i-1})}{c^2(U_i,U_{i-1})};\end{eqnarray*} \begin{eqnarray*}
a_{2,2}=\sum\limits_{i=1}^n\dfrac{4S^2_1(U_i)S^2_1(U_{i-1})}{c^2(U_i,U_{i-1})}&;& a_{2,3}=\sum\limits_{i=1}^n\dfrac{4C_2(U_{i})S_1(U_{i})C_2(U_{i-1})S_1(U_{i-1})}{c^2(U_i,U_{i-1})};\end{eqnarray*} \begin{eqnarray*}
a_{2,4}=\sum\limits_{i=1}^n\dfrac{4S_2(U_{i})S_1(U_{i})S_2(U_{i-1})S_1(U_{i-1})}{c^2(U_i,U_{i-1})}&;&
a_{3,3}= \sum\limits_{i=1}^n\dfrac{4C^2_2(U_i)C^2_2(U_{i-1})}{c^2(U_i,U_{i-1})};\\
a_{3,4}=\sum\limits_{i=1}^n\dfrac{4C_2(U_i)S_2(U_{i})C_2(U_{i-1})S_2(U_{i-1})}{c^2(U_i,U_{i-1})}\end{eqnarray*} \begin{eqnarray*}
a_{4,4}=\sum\limits_{i=1}^n\dfrac{4S_2(U_{i})S_1(U_{i})S_2(U_{i-1})S_1(U_{i-1})}{c^2(U_i,U_{i-1})}, \text{ $C_k(x)=\cos(2\pi kx),~ S_k(x)=\sin(2\pi k x)$}.\end{eqnarray*}
   
3. For the Legendre copula:
$\displaystyle
I_n(\Lambda)=\left(a_{j,k}
    \right)_{1\leq j,k\leq 2}
$ where \[a_{1,1}=\sum\limits_{i=1}^n\dfrac{9L^2(i)L^2(i-1)}{c^2(U_i,U_{i-1})};~ a_{2,2}=\sum\limits_{i=1}^n\dfrac{25P^2(i)P^2(i-1)}{c^2(U_i,U_{i-1})} \quad \text{and}\]  \[a_{1,2}=\sum\limits_{i=1}^n\dfrac{15
L(i)L(i-1)P(i)P(i-1)}{c^2(U_i,U_{i-1})}, \quad \text{ 
   $ L(i)=2U_i-1; P(i)=6U_i^2-6U_i+1$.}\]
\end{remark}

\subsection{The robust estimator of Longla and Peligrad (2021)}
~ The robust estimator defined in Longla and Peligrad (2021) uses an independent random sample to reduce the effect of the dependence. It makes use of kernel estimators to derive a central limit theorem that doesn't require knowledge of the variance of partial sums. For each of the three considered Markov chains, we consider the sequence $Y_k=(Y_{k1},\cdots,Y_{kn})$, where \begin{equation}\label{NewMarkov}
    Y_{ki}=\varphi_k(U_i)\varphi_k(U_{i-1}),\quad \mu_{Y_{k}}=\mathbb{E}Y_{ki}=\lambda_k.
\end{equation}
This sequence is a function of the stationary Markov chain $(\xi_0, \cdots, \xi_{n-1})$, with $\xi_i=(U_{i}, U_{i+1})$. 
According to Longla and Peligrad (2021) and considering the Gaussian kernel,  the robust estimator of $\mu_{Y_k}$, is 
\begin{equation}\label{estim3}
    \tilde{\lambda}_k=\dfrac{1}{nh_n}\sum\limits_{i=1}^n Y_{ki}\exp\left(-0.5(\dfrac{X_i}{h_n})^2\right),
\end{equation} 

\noindent where $X_i\sim\mathcal{N}(0,1)$ is the independently generated random sample,  the bandwidths sequence is   $h_n=\left[
        \dfrac{\overline{Y_{k}^2}}
{\overline{{Y}_{k}}^2n\sqrt{2}}
        \right]^{1/5}$ and estimators of means of $Y_{k,i}$ and $Y^2_{k,i}$ are given by $\overline{{Y}_{k}}=\dfrac{1}{n}\sum\limits_{i=1}^nY_{k,i}$
~and~ $\overline{{Y}_{k}^2}=\dfrac{1}{n}\sum\limits_{i=1}^n Y_{k,i}^2$, 
leading to the $(1-\alpha)100\%$ confidence intervals 
\begin{eqnarray}\label{ci3}
    \left(\tilde{\lambda}_k\sqrt{1+h_n^2}-z_{\alpha/2}\left(\dfrac{\overline{{Y}_{k}^2}}{nh_n\sqrt{2}}\right)^{1/2}, ~ \tilde{\lambda}_k\sqrt{1+h_n^2}+z_{\alpha/2}\left(\dfrac{\overline{{Y}_{k}^2}}{nh_n\sqrt{2}}\right)^{1/2}\right).
\end{eqnarray}

\begin{remark}
Given  that the means of $Y_{k,i}$ and $Y^2_{k,i}$ can also be computed as follows
\[\hat{\mathbb{E}}[Y_{k,i}]=\hat{\lambda}_k \text{ and } \hat{\mathbb{E}}[Y^2_{k,i}]=\hat{\mathbb{E}}(\varphi^2_k(U_0)\varphi^2_k(U_1))=1+\sum_{z=1}^s\hat{\lambda}_z(\int_0^1\varphi_z(x)\varphi_k^2(x) dx)^2,\]
the bandwidths sequence is provided by
      $\hat{h}_{n,k}=\left[
        \dfrac{\hat{\mathbb{E}}[Y^2_{k,i}]}{ \hat{\mathbb{E}}^2[Y_{k,i}]n\sqrt{2}}
\right]^{1/5}$, and the estimator of $\lambda_k$ becomes
\begin{equation}\label{estim3}
    \bar{\lambda}_k=\dfrac{1}{n\hat{h}_{n,k}}\sum\limits_{i=1}^n Y_{k,i}\exp\left(-0.5(\dfrac{X_i}{\hat{h}_{n,k}})^2\right)\end{equation}
which leads to the following  $(1-\alpha)100\%$ confidence interval
\begin{equation}\label{alt}
 \left(\bar{\lambda}_k\sqrt{1+\hat{h}_{n,k}^2}-z_{\alpha/2}\left(\dfrac{\hat{\mathbb{E}}[Y^2_{k,i}]}{n\hat{h}_{n,k}\sqrt{2}}\right)^{1/2}, ~ \bar{\lambda}_k\sqrt{1+\hat{h}_{n,k}^2}+z_{\alpha/2}\left(\dfrac{\hat{\mathbb{E}}[Y^2_{k,i}]}{n\hat{h}_{n,k}\sqrt{2}}\right)^{1/2}\right) 
\end{equation} where
\begin{enumerate}
    \item For the sine copula\\
    $\hat{\mathbb{E}}[Y^2_{1,i}]=1+\frac{\hat{\lambda}_2}{2}$,\quad  $\hat{\mathbb{E}}[Y^2_{2,i}]=1$,\quad  $\hat{h}_{n,1}=\left[\frac{2+\hat{\lambda}_2}{2n\sqrt{2}\hat{\lambda}_1^2}\right]^{1/5}$ and $\hat{h}_{n,2}=\left[\frac{1}{n\sqrt{2}\hat{\lambda}_2^2}\right]^{1/5}$.
\item For the sine-cosine copula\\
$\hat{\mathbb{E}}[Y^2_{1,i}]=1+\frac{\hat{\lambda}_2}{2}$,  $\hat{\mathbb{E}}[Y^2_{2,i}]=1$, $\hat{\mathbb{E}}[Y^2_{3,i}]=1+\frac{\hat{\lambda}_2}{2}$,  $\hat{h}_{n,4}=\left[\frac{1}{n\sqrt{2}\hat{\mu}_2^2}\right]^{1/5}$,
   
\noindent $\hat{\mathbb{E}}[Y^2_{4,i}]=1$, $\hat{h}_{n,1}=\left[\frac{2+\hat{\lambda}_2}{2n\sqrt{2}\hat{\lambda}_1^2}\right]^{1/5}$, $\hat{h}_{n,2}=\left[\frac{1}{n\sqrt{2}\hat{\lambda}_2^2}\right]^{1/5}$,
$\hat{h}_{n,3}=\left[\frac{2+\hat{\lambda}_2}{2n\sqrt{2}\hat{\mu}_1^2}\right]^{1/5}$.
\item For the Legendre copula\\
$\hat{\mathbb{E}}[Y^2_{1,i}]=1+\frac{4}{5}\hat{\lambda}_2$,  $\hat{\mathbb{E}}[Y^2_{2,i}]=1+\frac{20}{49}\hat{\lambda}_2$,  $\hat{h}_{n,1}=\left[\frac{5+4\hat{\lambda}_2}{5n\sqrt{2}\hat{\lambda}_1^2}\right]^{1/5}$ and $\hat{h}_{n,2}=\left[\frac{49+20\hat{\lambda}_2}{49n\sqrt{2}\hat{\lambda}_1^2}\right]^{1/5}$.
\end{enumerate}
    \end{remark}

\subsection {Data simulation and analysis of the examples}
We consider copulas with the following description for generation of Markov chains using the given algorithm. 
\begin{enumerate}
\item {\it  For the sine copula,  $\lambda_1=0.28$ and $\lambda_2=-0.15$;}
\item {\it For the sine-cosine copula, $\lambda_1=0.14, ~\lambda_2=0.13, ~\mu_1=-0.11, ~\mu_2=-0.12$;}
\item {\it For the Legendre copula, $\lambda_1=0.15$ and $\lambda_2=0.1$.}
\end{enumerate}
The obtained data is then used to estimate parameters and construct confidence intervals for the three cases studied above. We also provide comparisons with the MLE, the estimators $\hat{\Lambda}$ and $\tilde{\Lambda}$.
The obtained  confidence intervals are based on formulas (\ref{ci1}), (\ref{ci2}) and (\ref{ci3}). For the MLE, we have implemented an algorithm on $R$ using the level of significancy $\alpha=.05$. The graphs are obtained using Mathlab.

\subsubsection {Graphs of the log-likelihood functions}

\begin{figure}[h]
\begin{minipage}[c]{.46\linewidth}
     \begin{center}
             \includegraphics[width=5.5cm, height=4cm]{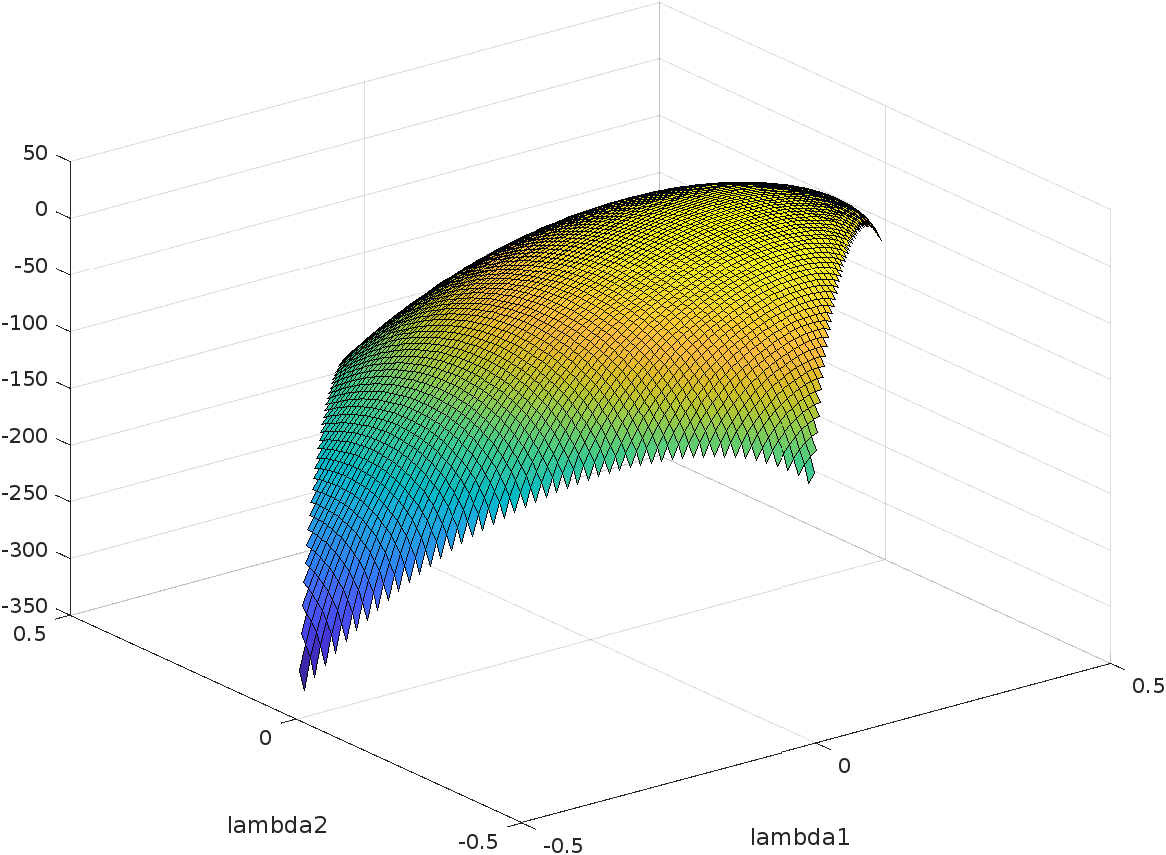}
         \end{center}
   \end{minipage} \hfill
   \begin{minipage}[c]{.46\linewidth}
    \begin{center}
            \includegraphics[width=5.5cm, height=4cm]{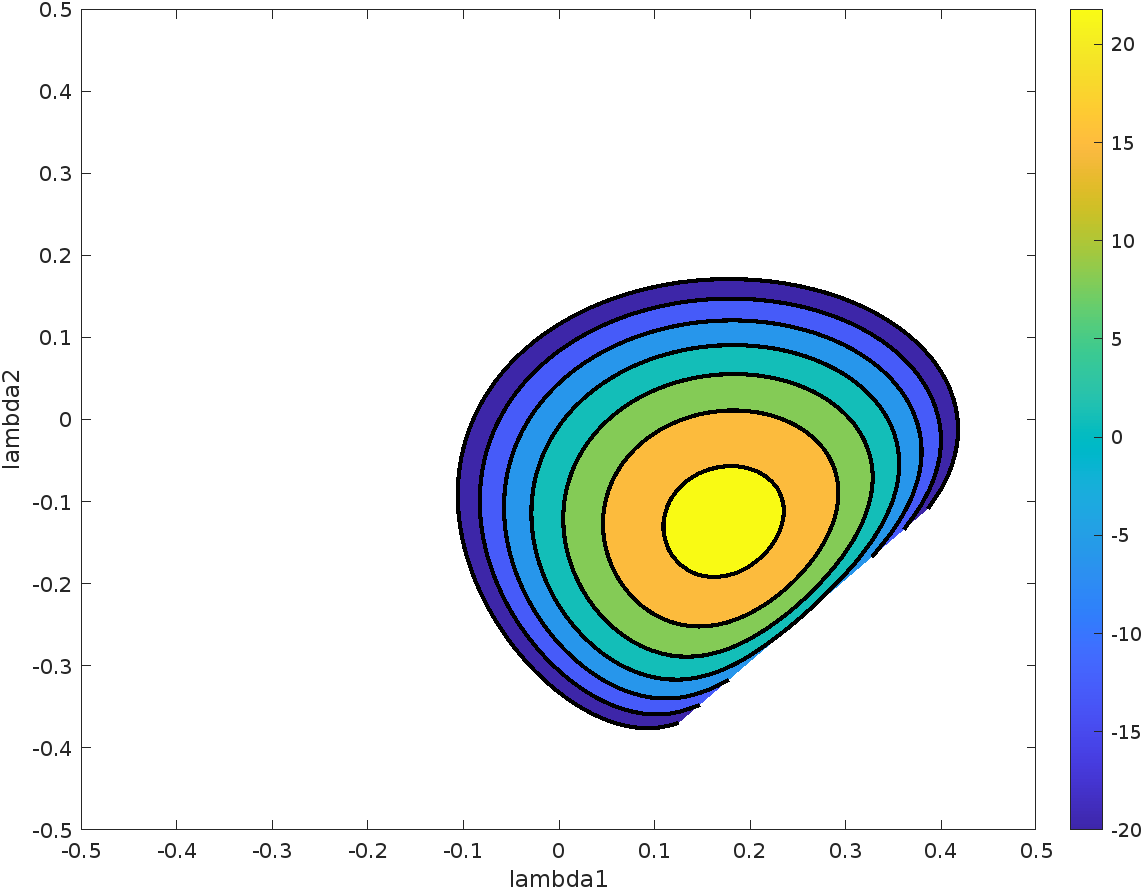}
        \end{center}
 \end{minipage}
 \caption{Log-likelihood function under thesine copula}\label{f1}
 \end{figure}

\begin{figure}[h]
\begin{minipage}[c]{.46\linewidth}
     \begin{center}
             \includegraphics[width=5.5cm, height=4cm]{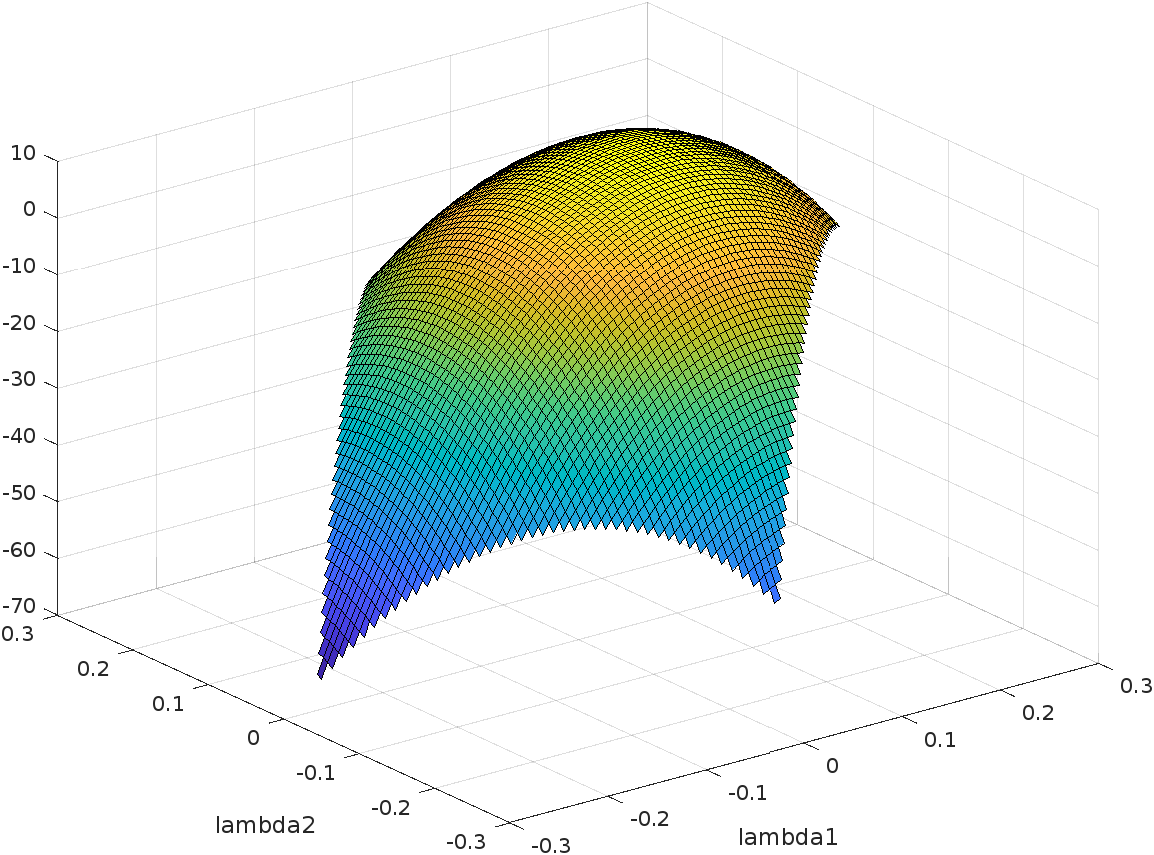}
         \end{center}
   \end{minipage} \hfill
   \begin{minipage}[c]{.46\linewidth}
    \begin{center}
            \includegraphics[width=5.5cm, height=4cm]{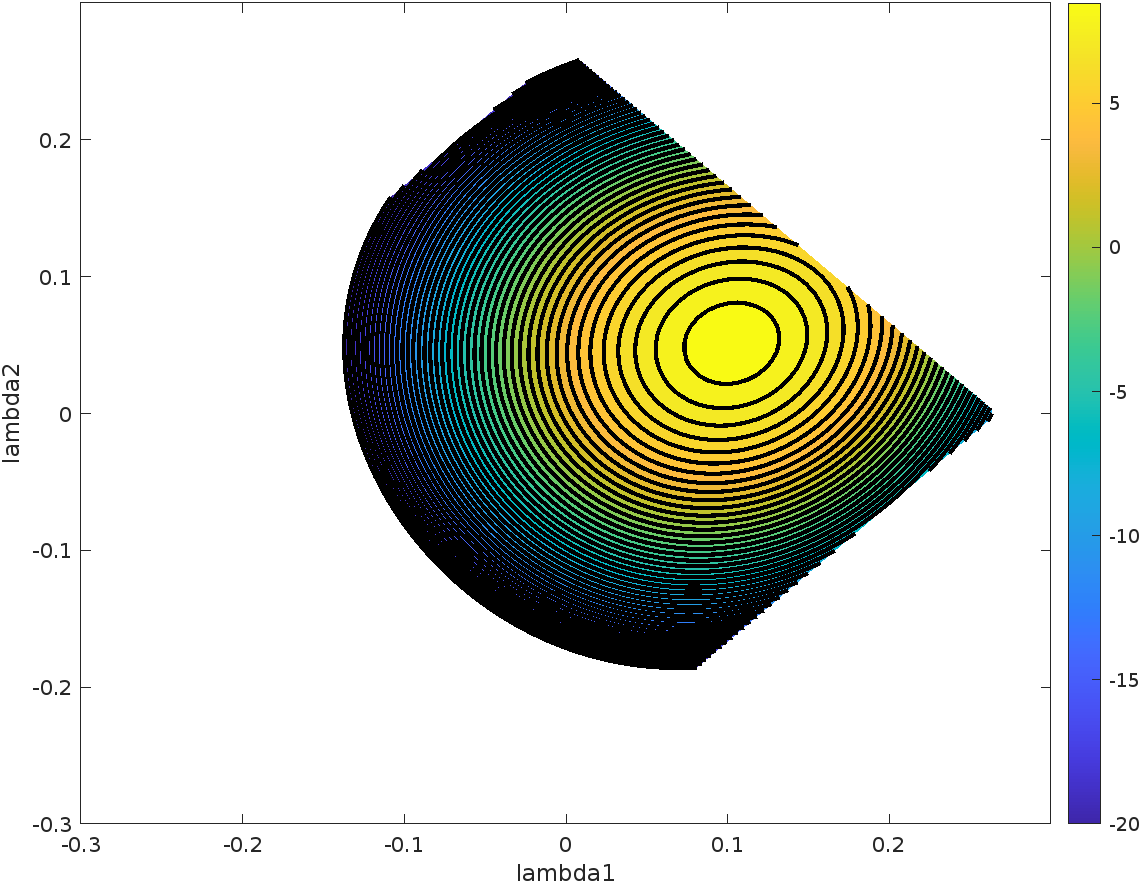}
        \end{center}
 \end{minipage}
 \caption{Log-likelihood function under the sine-cosine copula}\label{f2}
 \end{figure}

\begin{figure}[h]
\begin{minipage}[c]{.46\linewidth}
     \begin{center}
             \includegraphics[width=5.5cm, height=4cm]{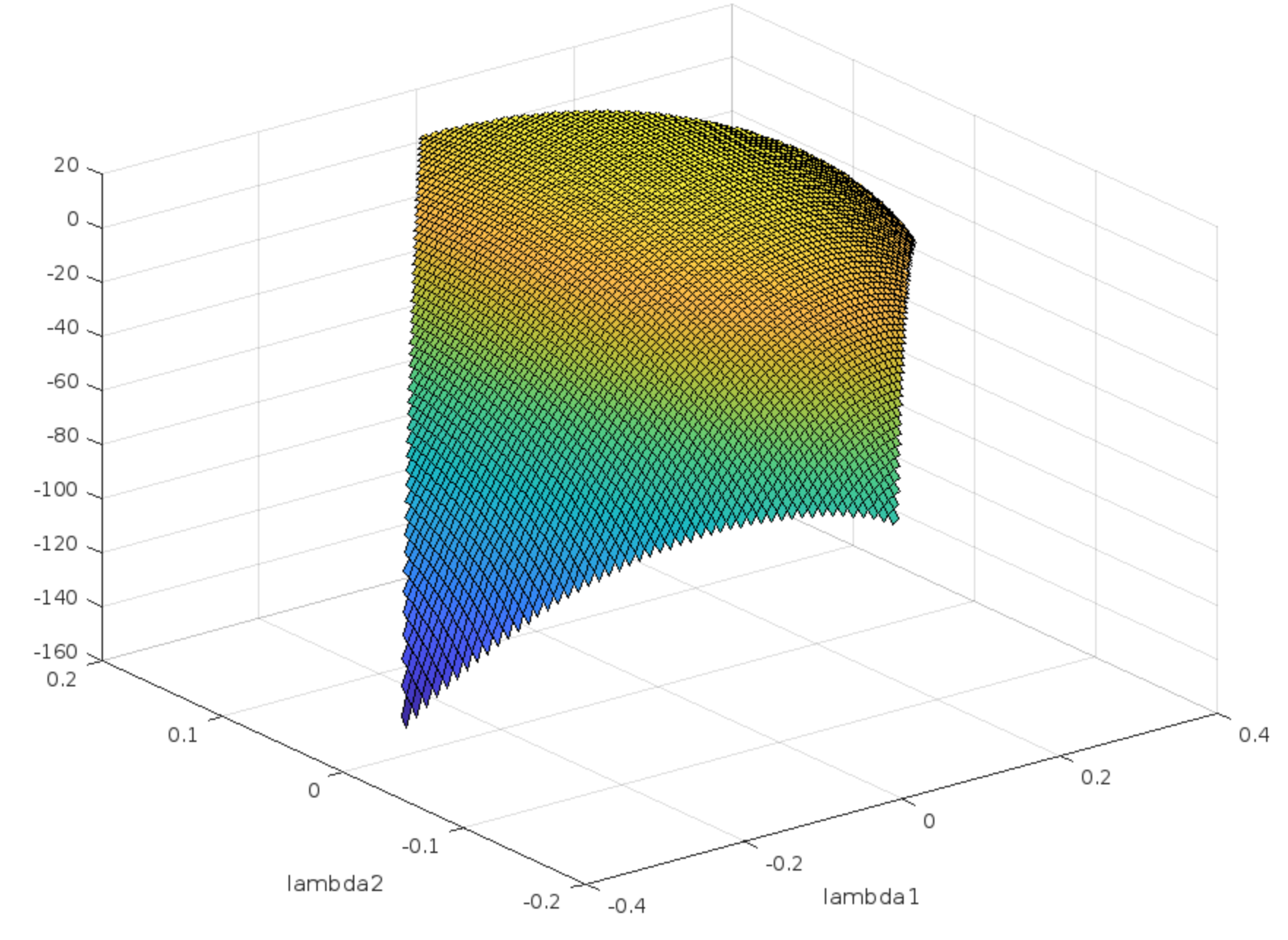}
         \end{center}
   \end{minipage} \hfill
   \begin{minipage}[c]{.46\linewidth}
    \begin{center}
            \includegraphics[width=5.5cm, height=4cm]{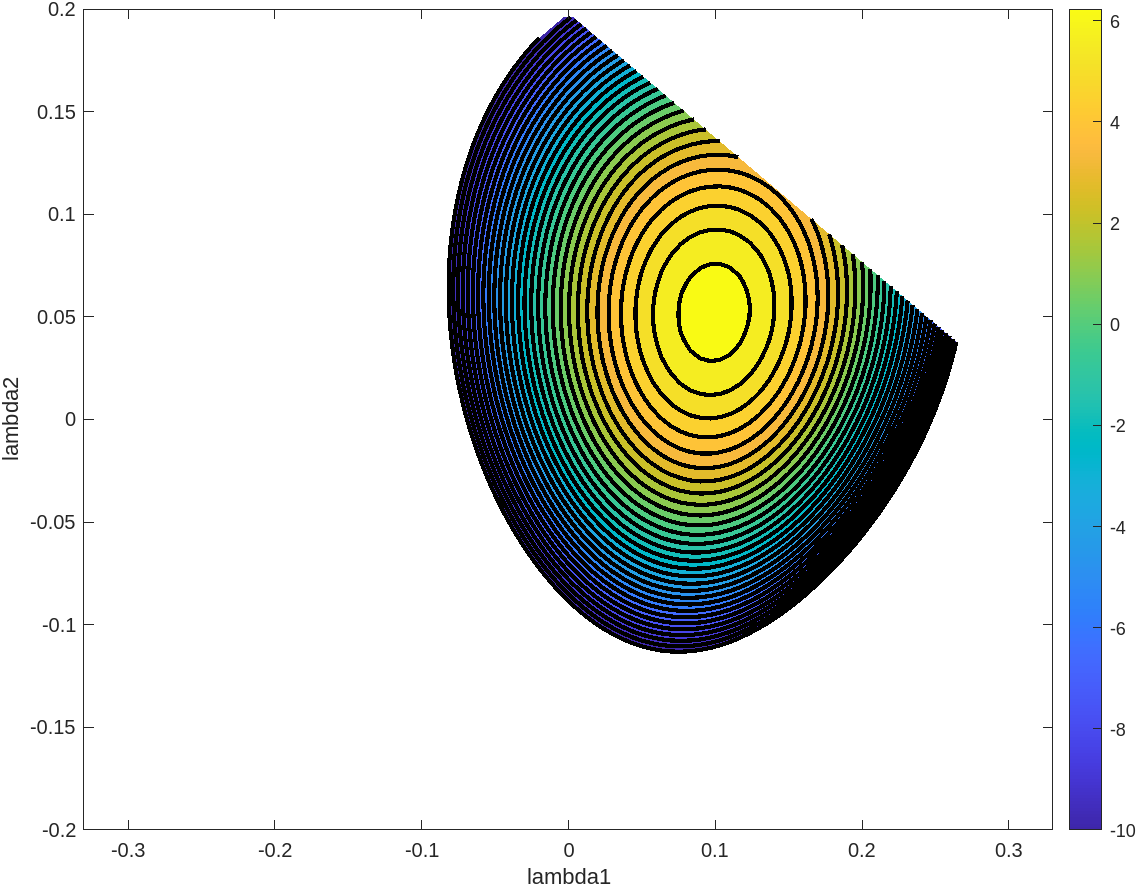}
        \end{center}
 \end{minipage}
 \caption{Log-likelihood function under the Legendre copula}\label{f3}
 \end{figure}

 \begin{remark}{Comments on the log-likelihood graphs}
\begin{enumerate}
\item   The log-likelihood function of the sine copula in figure \ref{f1} is obtained using the Markov generated with parameter $\Lambda=(0.14,-0.8)$;
\item   For the sine-cosine copula, we use the Markov chain generated with parameter $\Lambda=(0.07,0.065,-0.055,-0.06)$. In order to obtain the representation in figure \ref{f2}, we've fixed $\mu_1$ and $\mu_2$ in the likelihood function;
  \item The log-likelihood function of the Legendre copula shown in Figure \ref{f3} uses the Markov chain generated with parameter $\Lambda=(0.075, 0.05)$.\end{enumerate} 

In each of the cases, the contour plot of the log-likelihood function exhibits evident uniqueness of the MLE. This is shown by spectral colors of the level curves stretching towards a single interior point. This observation provides reassurance regarding reliability of the MLE and confirms the theoretical insights.
    \end{remark}


\subsubsection{Analysis of simulation results}
For every given sample size, a data set has been simulated, representing a Markov chain generated by the considered copula and the uniform marginal distribution. Computations have been made for each of the copulas in $R$. 
After running 100 times each of these simulations, we have counted the number of intervals including the true value of the parameter to get the coverage probability, wich should represent  $95$ out of the $100$ constructed intervals per case. We also calculated the mean length of the  confidence interval. Conclusions of the simulation study are given in the remark below. In these tables, we use CP for coverage probability and CIML for confidence interval's mean length.

\begin{table}[H]
\centering
\begin{tabular}[t]{|c|c|c|c|c|c|c|}
\hline
 \multicolumn{7}{|c|}{ Sine-cosine copula for  $\lambda_1=0.14,~\lambda_2=0.13,~ \mu_1=-0.11,~\mu_2=-0.12$}
 \\
\hline
\multirow{2}{*}{Estimator}&\multicolumn{2}{|c|}{n=4999}&\multicolumn{2}{|c|}{n=9999}&\multicolumn{2}{|c|}{n=19999}\\
 &CP&CIML&CP&CIML&CP&CIML\\
\hline\hline
{$\hat{\lambda}_1=0.1504$}
&97&0.0523&97&0.0369&95&0.0256 \\
\hline
$\hat{\lambda}_1^{ML}=0.1349$
&98&0.0543&98&0.0361&98&0.0255\\
\hline
$\tilde{\lambda}_1=0.1278$
&98&0.0517&93&0.0525&95&0.0407\\
\hline
$\bar{\lambda}_1=0.1280$
&96&0.0549&95&0.0560&93&0.0415\\
\hline\hline
{$\hat{\lambda}_2=0.1339$}
&98&0.0511&97&0.0380&98&0.0272 \\
\hline
$\hat{\lambda}_2^{ML}=0.1315$
&95&0.0497 &98&0.0362&97&0.0253\\
\hline
$\tilde{\lambda}_2$=0.1281
&94&0.0509&92&0.0515&94&0.0398\\
\hline
$\bar{\lambda}_2$=0.1130
&90&0.0529&0.575&99&94&0.0413\\
\hline\hline
{$\hat{\mu_1}=-0.1039$}
&96&0.0732&91&0.0347&93&0.0250 \\
\hline
${\hat{\mu}_1^{ML}}=-0.1023$
&95&0.0694&95&0.0360&93&0.0249\\
\hline
$\tilde{\mu}_1$=0.0945
&92&0.0665&92&0.0506&92&0.0384\\
\hline
$\bar{\mu}_1=-0.1067$
&93&0.0702 &96&0.0522&96&0.0396\\
\hline\hline
{$\hat{\mu_2}=-0.1223$}
&99&0.0744&94&0.0368&90 &0.0249\\
\hline
$\hat{\mu}_2^{ML}=-0.1233$
&99&0.0683&95&0.0359&93 &0.0248\\
\hline
$\tilde{\mu}_2=-0.1042$
&95&0.0662&95&0.0530&93&0.0395\\
\hline
 $\bar{\mu}_2=-0.1294$
&93&0.0678&93&0.0513&93&0.0392\\
\hline
\end{tabular}
\caption{Estimates, coverage probabilities and mean lengths.}
\label{tt}
\end{table}
\begin{remark}[On simulations]{}{}
The following summarizes the results of the conducted simulations.
    \begin{enumerate}
    \item In Table \ref{tt}, the estimates are provided for a sample size of $10000$. For Table \ref{cosinetab}, we utilize a sample size of $5000$, while for Table \ref{legendretab}, a sample size of $20000$ is considered for the estimates.
    \item Through $R$'s {\it constrOptim} function, the MLE is obtained. This function is initialized at $(0,0)$ for the sine and Legendre copulas and at  $(0,0,0,0)$ for the sine-cosine copula. We successfully optimized the likelihood function while maintaining the linear inequality constraints thanks to the function's flexibility.
        \item For each type of copula, the proposed estimator $\hat{\Lambda}$ challenges the MLE by demonstrating its effectiveness in providing reliable parameter estimates in situations involving small sample sizes. Furthermore, it exhibits an increased level of precision as the sample size grows larger.

\begin{table}[ht]
\centering
 \begin{tabular}{|c|c|c|c|c|c|c|}
\hline
 \multicolumn{7}{|c|}{ Legendre copula for $ \lambda_1=0.15,~\lambda_2=0.1$}
 \\
\hline
\multirow{2}{*}{Estimator}&\multicolumn{2}{|c|}{n=4999}&\multicolumn{2}{|c|}{n=9999}&\multicolumn{2}{|c|}{n=19999}\\
 &CP&CIML&CP&CIML&CP&CIML\\
\hline
\hline
 $\hat{\lambda}_1$=0.1511
&93&0.0522&94&0.0387&94&0.0274 \\
\hline
$\hat{\lambda}_1^{ML}$=0.1512
&96&0.0535&97&0.0383&96&0.0268\\
\hline
$\tilde{\lambda}_1$=0.1309
&92&0.0737&97&0.0586&91&0.0418\\
\hline
$\bar{\lambda}_1$= 0.1555
&92&0.0733&95&0.0572&93&0.0426\\
\hline\hline
$\hat{\lambda}_2$=0.1101
&88&0.0502&89&0.0359&92&0.0234\\
\hline
 $\hat{\lambda}_2^{ML}$=0.1081
&99&0.0539&99&0.0381&91&0.0248\\
\hline
 $\tilde{\lambda}_2$=0.1042
&90&0.0665&91&0.0509&88&0.0373\\
\hline
$\bar{\lambda}_2$=0.1042
&90&0.0600&91&0.0476&93&0.0375\\
\hline
\end{tabular}
\caption{Estimates, coverage probabilities and mean lengths.}
\label{cosinetab}
\end{table}

        \item One significant advantage of $\hat{\Lambda}$ over $\Lambda^{ML}$ is that the latter necessitates more computational time due to the absence of a closed of Fisher information matrix  for the MLE and simple form of the variance for and estimators of our proposed method for considered copulas. For each of the cases $\hat{\Lambda}$ has an explicit expression for both the estimate and the confidence interval, resulting in reduced computation time  and reduced errors for the algorithm.

\begin{table}[ht]
\centering
\begin{tabular}{|c|c|c|c|c|c|c|}
\hline
 \multicolumn{7}{|c|}{ Sine copula for $  \lambda_1=0.28,~\lambda_2=-0.15$}
 \\
\hline
\multirow{2}{*}{Estimator}&\multicolumn{2}{|c|}{n=4999}&\multicolumn{2}{|c|}{n=9999}&\multicolumn{2}{|c|}{n=19999}\\
 &CP&CIML&CP&CIML&CP&CIML\\
\hline\hline
$\hat{\lambda}_1$=0.2998
&94&0.0417&98&0.0241&97&0.0127\\
\hline
$\hat{\lambda}_1^{ML}=
0.2991$
&95&0.0298&96&0.0084&97&0.0013\\
\hline
$\tilde{\lambda}_1=0.2784$
&91&0.0662&97&0.0455&91&0.0290\\
\hline$\bar{\lambda}_1=
0.3078$
&97&0.0699&98&0.0491&97&0.0334\\
\hline\hline
$\hat{\lambda}_2$ =$-0.1354$
&97&0.0235&95&0.0058&96&0.0003\\
\hline
$\hat{\lambda}_2^{ML}=0.1438$ &97&0.0044&94&0.0003&95 &0.0002\\
\hline
$\tilde{\lambda}_2=-0.1213$
&98&0.0356&94&0.0202&96&0.0071\\
\hline
$\bar{\lambda}_2=-0.1320$
&97 &0.0383&93&0.0236&95&0.0053\\
\hline
\end{tabular}
\caption{Estimates, coverage probabilities and mean lengths.}
\label{legendretab}
\end{table}

        \item In Table \ref{legendretab}, the coverage probabilities of the MLE are unexpectedly high for small sample sizes, given the specified significance level. This can be interpreted as an instance of over-coverage, indicating that the MLE provides intervals that are wider than necessary to achieve the desired level of confidence.
        \item The alternative estimator $\bar{\Lambda}$, in comparison to $\tilde{\Lambda}$, offers improved point estimates and coverage probabilities for each copula. It demonstrates enhanced accuracy in estimating the true parameter values and achieves more reliable coverage probabilities within the considered framework.
    \end{enumerate}
\end{remark}

\section{Conclusions and remarks}
This paper studies copula-based Markov chains with uniform marginals. These Markov chains are based, in general, on copulas of Longla (2023). This is the first study of properties of estimators of parameters of these copulas. A multivariate central limit theorem is provided for copula parameters estimators. These copulas have densities of the general form \eqref{cop}, which eases the investigation of asymptotic properties via mixing. This central limit theorem has been applied to examples of sine copulas, sine-cosine copulas and Legendre copulas introduced in Longla (2023). The proposed estimators of parameters have been compared to the robust estimators of Longla and Peligrad (2021). This theorem has been used, to show that our estimators perform better than the MLE for sample sizes above $5000$.

The simulation study that has been performed using $R$, shows that it takes more time to run the procedures of MLE than to run their equivalents for our proposed estimators. Moreover, the possibility of having several points that indicate the maximum of the likelihood functions, and the difficulty to obtain the true variance of the MLE, are indicators that our estimators are better tools.  Note that it takes long iterations to approximate integrals that lead to the variance-covariance matrix of the asymptotic distribution of the MLE, while the variance-covariance matrix of our estimator is know in closed form. It is also good to note that on large samples, coverage probabilities also indicate that the MLE is not better. 

One of the main takeaways of this work is the provided theoretical proof of large sample properties of our estimators and the MLE for Markov chains generated by the copulas with densities of the form \eqref{cop}. This opens a path to several other questions that can be investigated for these copula families. Future research on the topic includes the use of these theorems to investigates tests of equality of two or more bivariate distributions based on realizations of Markov chains that each of them generates. We also consider working on general estimation problems, when the functions $\varphi_k(x)$ have a parameter. All these issues can also be included in problems with a general marginal distribution. The marginal distribution can be considered with parametric or non-parametetric form. All these questions are interesting and not obvious, and are among our topics of current research.

Another extension that can be considered here is the work on such questions for a perturbation of a copula other than the independence copula. The theoretical  work on this last topic is a little more complicated, because orthogonality becomes more complex.

\section*{Statements and declarations} 
All authors certify that they have no affiliations with or involvement in any organization or entity with any financial interest or non-financial interest in the subject matter or materials discussed in this manuscript. They have no competing interest and have no funding to disclose.

\newpage
\section{Appendix of proofs of major results}

\subsection{Proof of Theorem \ref{lambda1}}
Assume $U_0, U_1, \cdots, U_n $ is a copula-based Markov chain generated by a copula of the form \eqref{cop}. The stationary distribution of the Markov chain is uniform on $(0,1)$. By construction, $\mathbb{E}(\varphi_j(U_i))=0$, $\mathbb{E}(\varphi_j(U_i)\varphi_l(U_{i}))=\delta_{il}$ for any $U_i$. Here $\delta_{il}=1$ when $i=l$ and $\delta_{il}=0$ when $i\ne l$. It follows that $\mathbb{E}(\hat{\lambda}_k)=\lambda_k$.
Moreover, $$var(\hat{\lambda}_k)=\frac{1}{n^2}\sum_{i=1}^n\sum_{j=1}^ncov(\varphi_k(U_i)\varphi_k(U_{i-1}), \varphi_k(U_j)\varphi_k(U_{j-1}))$$ $$=\frac{1}{n^2}\sum_{i=1}^n\sum_{j=1}^n\mathbb{E}(\varphi_k(U_i)\varphi_k(U_{i-1}) \varphi_k(U_j)\varphi_k(U_{j-1}))-\lambda_k^2.$$
Using stationarity and reversibility, we obtain $$var(\hat{\lambda}_k)=\frac{1}{n}\mathbb{E}(\varphi^2_k(U_1)\varphi^2_k(U_{0}))+\frac{2(n-1)}{n^2}\mathbb{E}(\varphi_k(U_{0}) \varphi^2_k(U_1)\varphi_k(U_2))+$$$$+\frac{1}{n^2}\sum_{s=2}^{n-1}2(n-s)\mathbb{E}(\varphi_k(U_0)\varphi_k(U_{1}) \varphi_k(U_s)\varphi_k(U_{s+1}))=\lambda_k^2.$$
To end these computations, we need the three joint distributions of $(U_0, U_1, U_2)$ and $(U_0,U_1,U_s, U_{s+1})$ for $s\ge 2$.  The following Lemma can be easily established using the definitions.
\begin{lemma} \label{comput}
If $(U_0,U_1,\cdots, U_n)$ is a copula-based Markov chain generated by the copula with density \eqref{cop}, then the following holds: 

\begin{enumerate}
\item  The density of $(U_0,U_1,U_2)$ is $\displaystyle c(u_0, u_1,u_2)=c(u_0,u_1)c(u_1,u_2)=$ $$=1+\sum\lambda_m \varphi_m(u_0)\varphi_m(u_1)+\sum\lambda_n \varphi_z(u_1)\varphi_z(u_2)$$ $$+ \sum\sum\lambda_m\lambda_z \varphi_m(u_0)\varphi_m(u_1)\varphi_z(u_1)\varphi_z(u_2).$$ 
\item  $\displaystyle \mathbb{E}(\varphi_k(U_{0}) \varphi^2_k(U_1)\varphi_k(U_2))=\lambda_k^2\int_{0}^1 \varphi^4_k(x)dx$. 
\item $\displaystyle \mathbb{E}(\varphi^2_k(U_1)\varphi^2_k(U_{0}))=1+\sum \lambda_z (\int_0^1\varphi_z(x)\varphi_k^2(x)dx)^2.$
\item The joint density of $ (U_0,U_1,U_s, U_{s+1})$ is $$\displaystyle c(u_0, u_1,u_s,u_{s+1})=c(u_0,u_1)(1+\sum \lambda_z^{s-1}\varphi_z(u_1)\varphi_z(u_s))c(u_s,u_{s+1}).$$
\item $\mathbb{E}(\varphi_k(U_0)\varphi_k(U_{1}) \varphi_k(U_s)\varphi_k(U_{s+1})\lvert{}U_0,U_1,U_s)=\lambda_k\varphi_k(U_0)\varphi_k(U_{1}) \varphi^2_k(U_s).$
 \item The joint density of $ (U_0,U_1,U_s)$ is $$\displaystyle c(u_0, u_1,u_s,u_{s+1})=c(u_0,u_1)(1+\sum \lambda_z^{s-1}\varphi_z(u_1)\varphi_z(u_s)).$$
\item $\mathbb{E}(\varphi_k(U_0)\varphi_k(U_{1}) \varphi_k(U_s)\varphi_k(U_{s+1}))=\lambda_k\mathbb{E}(\varphi_k(U_0)\varphi_k(U_{1}) \varphi^2_k(U_s))$ $$ =\lambda_k^2\mathbb{E}(\varphi^2_k(U_{1}) \varphi^2_k(U_s)).$$
\item $\mathbb{E}(\varphi_k(U_0)\varphi_k(U_{1}) \varphi_k(U_s)\varphi_k(U_{s+1}))=\lambda_k^2(1+\sum\lambda^{s-1}_z(\int_0^1\varphi_z(x)\varphi_k^2(x)dx)^2).$
\end{enumerate}
\end{lemma}
 Based on Lemma \ref{comput}, we can conclude that 
$$var(\hat{\lambda}_k)=\frac{1}{n}(1-\lambda_k^2+\sum \lambda_z (\int_0^1\varphi_z(x)\varphi_k^2(x)dx)^2)+\frac{2(n-1)\lambda_k^2}{n^2}(\int_0^1\varphi^4_k(x)dx-1)$$
\begin{equation}
+\sum_{s=2}^{n-1}\frac{2(n-s)\lambda_k^2}{n^2}(\sum\lambda^{s-1}_z(\int_0^1\varphi_z(x)\varphi_k^2(x)dx)^2).
\end{equation}
Therefore, knowing that $\sum \lvert{}\lambda_z\lvert{} (\int_0^1\varphi_z(x)\varphi_k^2(x)dx)^2)=K<\infty$ and $\int_{0}^1\varphi_k^4(x)dx<\infty$, we have
\begin{equation}
\lim_{n\to\infty}var(\hat{\lambda}_k)=\lim_{n\to\infty}\frac{2\lambda_k^2}{n^2}\sum_{s=3}^{n-1}(n-s)(\sum\lambda^{s-1}_z(\int_0^1\varphi_z(x)\varphi_k^2(x)dx)^2).
\end{equation}
Provided that $\underline{\lambda}=\sup_{z}\lvert{}\lambda_z\lvert{}<1$ , we obtain $$\lim_{n\to\infty}var(\hat{\lambda}_k)\leq \lim_{n\to\infty}\frac{2\lambda_k^2}{n^2}\sum_{s=3}^{n-1}(n-s)\underline{\lambda}^{s-2}(\sum\lvert{}\lambda_z\lvert{}(\int_0^1\varphi_z(x)\varphi_k^2(x)dx)^2), $$ $$\lim_{n\to\infty}var(\hat{\lambda}_k)\leq \lim_{n\to\infty}\frac{2K\lambda_k^2}{n^2}\sum_{s=3}^{n-1}(n-s)\underline{\lambda}^{s-2}.$$

Moreover, the last sum can be written as $$\sum_{s=3}^{n-1}(n-s)\underline{\lambda}^{s-2}=(n\sum_{s=3}^{n-1}\underline{\lambda}^{s-2} -\frac{1}{\underline{\lambda}}\frac{d}{d\underline{\lambda}}(\sum_{s=3}^{n-1}\underline{\lambda}^{s}))=$$ $$\frac{n\underline{\lambda}(1-\underline{\lambda}^{n-2})}{1-\underline{\lambda}}-\frac{1}{\underline{\lambda}}\frac{d}{d\underline{\lambda}}\frac{\underline{\lambda}^2(1-\underline{\lambda}^{n-2})}{1-\underline{\lambda}}= \frac{n(\underline{\lambda}-\underline{\lambda}^{n-1})-(2-n\underline{\lambda}^{n-2})}{1-\underline{\lambda}}+\frac{\underline{\lambda}^2-\underline{\lambda}^n}{(1-\underline{\lambda})^2}. $$

If $\sup_z \lvert{}\lambda_z\lvert{}<1$, then as $n>2 $ this sum behaves like $\displaystyle n\frac{\underline{\lambda}}{1-\underline{\lambda}}.$ 
Therefore $\lim_{n\to\infty}var(\hat{\lambda}_k)=0$. This implies, that under the conditions of the statement, $\hat{\lambda}_k$ is a consistent estimator of $\lambda_k$. Moroever, from this proof, it follows that $n var(\hat{\lambda}_k)\to \sigma^2$ with $0< \sigma^2< \infty$. We can also establish easily that 
$Y_i=(X_{i-1}, X_i), i=1,\cdots,n $ forms a reversible Markov chain. The variables $Z_i=\varphi_k(U_{i})\varphi_k(U_{i-1}):=f(Y_i)$ satisfy $n^{-1} var(\sum Z_i)\to \sigma^2$ because $\sum Z_i=n\hat{\lambda}_k$. Therefore, the Kipnis and Varadhan (1986) central limit theorem holds for $\hat{\lambda}_k$.

\subsection{Proof of Theorem \ref{lambdam}}
Assume $(U_0,U_1, \cdots, U_n)$ is a copula-based Markov chain generated by a copula of the form \eqref{cop}.
It is already clear that each of the $\hat{\lambda}_k$ is a consistent estimator of $\lambda_k$. To prove this theorem, all we need to do now is to show by the Cram\'er-Wold device that for every ${\bf t}=(t_1, \cdots, t_s)'$, $\sqrt{n}({\bf t}'\hat{\Lambda}-{\bf t}'\Lambda)\to N({\bf 0}, {\bf t}'\Sigma {\bf t})$. And by the Kipnis and Varadhan central limit theorem for reversible Markov chains, to complete this proof, we need to show that $n var(t'\hat{\lambda})\to t'\Sigma t$ as $n\to\infty$. But it is easy to realize that 
$$n var({\bf t}'\hat{\Lambda})=\sum_{j=1}^s\sum_{i=1}^{s}t_jt_i n cov(\hat{\lambda}_{k_i}, \hat{\lambda}_{k_j}).$$
Therefore, to prove the central limit theorem, it is enough to show that $n cov(\hat{\lambda}_{k_i}, \hat{\lambda}_{k_j})\to \Sigma_{ij}$ as $n\to\infty$. As we have already shown this convergence for $i=j$, it is enough to conduct the proof for $i\ne j$.
But, $$n cov(\hat{\lambda}_{k_i}, \hat{\lambda}_{k_j})=\frac{1}{n}\sum_{p=1}^{n}\sum_{q=1}^{n}cov(\varphi_{k_i}(U_{p-1})\varphi_{k_i}(U_{p}), \varphi_{k_j}(U_{q-1})\varphi_{k_j}(U_{q})).$$

To compute this covariance, we consider $4$ separate cases $(p=q, p=q-1, q=p-1, \lvert{}p-q\lvert{}=m>1 )$. The cases $p=q-1$ and $q=p-1$ are equivalent due to symmetry. 
\begin{enumerate}
\item For $p=q$ and $i\ne j$, using stationarity we have $$cov(\varphi_{k_i}(U_{p-1})\varphi_{k_i}(U_{p}), \varphi_{k_j}(U_{q-1})\varphi_{k_j}(U_{q}))=$$ $$\mathbb{E}(\varphi_{k_i}(U_{0})\varphi_{k_i}(U_{1})\varphi_{k_j}(U_{0})\varphi_{k_j}(U_{1}))-\lambda_{k_j}\lambda_{k_i}=$$
$$=\sum_{z=1}^{\infty}\lambda_z (\int_0^1\varphi_{z}(x)\varphi_{k_i}(x)\varphi_{k_j}(x)dx)^2-\lambda_{k_j}\lambda_{k_i}.$$
\item For $q=p-1$ and $i\ne j$, using stationarity we have 
$$cov(\varphi_{k_i}(U_{p-1})\varphi_{k_i}(U_{p}), \varphi_{k_j}(U_{p})\varphi_{k_j}(U_{p+1}))=$$ $$\mathbb{E}(\varphi_{k_i}(U_{0})\varphi_{k_i}(U_{1})\varphi_{k_j}(U_{1})\varphi_{k_j}(U_{2}))-\lambda_{k_j}\lambda_{k_i}=$$ 
$$=\lambda_{k_j}\lambda_{k_i}(\int_0^1\varphi^2_{k_i}(x)\varphi^2_{k_j}(x)dx-1).$$
\item For $\lvert{}p-q\lvert{}=m\ge 2$ and $j\ne j$, using stationary we have
$$cov(\varphi_{k_i}(U_{p-1})\varphi_{k_i}(U_{p}), \varphi_{k_j}(U_{p+m-1})\varphi_{k_j}(U_{p+m}))=$$ $$\mathbb{E}(\varphi_{k_i}(U_{0})\varphi_{k_i}(U_{1})\varphi_{k_j}(U_{m})\varphi_{k_j}(U_{m+1}))-\lambda_{k_j}\lambda_{k_i}=$$ $$\lambda_{k_j}\lambda_{k_i}(\mathbb{E}(\varphi^2_{k_i}(U_{1})\varphi^2_{k_j}(U_{m}))-1)=$$ $$ =\lambda_{k_j}\lambda_{k_i}\sum_{z=1}^{\infty}\lambda_z^{m-1} (\int_0^1\varphi_{z}(x)\varphi_{k_i}^2(x)dx)(\int_0^1\varphi_{z}(x)\varphi_{k_j}^2(x)dx).$$
\end{enumerate}

Similarly to Haggstrom and Rosenthal (2007), it follows that $$\lim_{n\to \infty} n cov(\hat{\lambda}_{k_i}, \hat{\lambda}_{k_j})=cov(\varphi_{k_i}(U_{0})\varphi_{k_i}(U_{1}), \varphi_{k_j}(U_{0})\varphi_{k_j}(U_{1}))+$$ $$2\sum_{m=1}^{\infty}cov(\varphi_{k_i}(U_{0})\varphi_{k_i}(U_{1}), \varphi_{k_j}(U_{m})\varphi_{k_j}(U_{m+1})).$$
Knowing that $c(u,v)$ is square integrable, we can establish easily that for $m\le 3$ each of the terms of the series is in absolute value smaller than $\underline{\lambda}^{m-3} K$, where $\underline{\lambda}$ is the suppremum of $\lvert{}\lambda_z\lvert{}$ over $z$ and $K$ is a constant that depends only on $k_i,k_j$. This implies that the series with sum starting at $m=3$ is convergent (using geometric series to end the argument as $0<\underline{\lambda}<1$). Moreover the term with $m=1$ is finite under any conditions and the term with $m=2$ is finite if and only if 
$\sum_{z=1}^{\infty}\lambda_z (\int_0^1\varphi_{z}(x)\varphi_{k_i}^2(x)dx)(\int_0^1\varphi_{z}(x)\varphi_{k_j}^2(x)dx)<\infty.$
On the other hand, $cov(\varphi_{k_i}(U_{0})\varphi_{k_i}(U_{1}), \varphi_{k_j}(U_{0})\varphi_{k_j}(U_{1}))$ is finite if and only if 
$$\sum_{z=1}^{\infty}\lambda_z (\int_0^1\varphi_{z}(x)\varphi_{k_i}(x)\varphi_{k_j}(x)dx)^2<\infty.$$

\subsection{Proof of Proposition \ref{prop1}}
Following the proof of Theorem \ref{lambdam} and Haggstorm amd Rosenthal (2007), we have
    \begin{eqnarray}\label{covv}
    \underset{k\neq j}{\Sigma_{k,j}}&=&\underset{n\rightarrow\infty}{\lim~}ncov(\hat{\lambda}_k,\hat{\lambda}_j)\nonumber\\
&=&cov(\varphi_k(U_0)\varphi_k(U_1),\varphi_j(U_0)\varphi_j(U_1))+2cov(\varphi_k(U_0)\varphi_k(U_1),\varphi_j(U_1)\varphi_j(U_2))\nonumber\\ &&+2\sum\limits_{m=2}^\infty cov(\varphi_k(U_0)\varphi_k(U_1),\varphi_k(U_m)\varphi_k(U_{m+1}).
\end{eqnarray}

Using the obtained expected values, we derive the following.
\begin{eqnarray*}    cov(\varphi_k(U_0)\varphi_k(U_1),\varphi_j(U_0)\varphi_j(U_1))&=&\mathbb{E}\left(\varphi_k(U_0)\varphi_k(U_1)\varphi_j(U_0)\varphi_j(U_1)\right)-\lambda_k\lambda_j\nonumber\\
&=&\sum\limits_{z=1}^s\lambda_z\left(\int_0^1\varphi_z(x)\varphi_k(x)\varphi_j(x)dx\right)^2-\lambda_k\lambda_j.\\
    cov(\varphi_k(U_0)\varphi_k(U_1),\varphi_j(U_1)\varphi_j(U_2))
&=&\mathbb{E}\left(\varphi_k(U_0)\varphi_k(U_1)\varphi_j(U_1)\varphi_j(U_2)\right)-\lambda_k\lambda_j\nonumber\\&=&\lambda_k\lambda_j\left(\int_0^1\varphi^2_k(x)\varphi_j^2(x)-1\right).
\end{eqnarray*}
\begin{eqnarray*}
cov(\varphi_k(U_0)\varphi_k(U_1),\varphi_k(U_m)\varphi_k(U_{m+1})=& {} \\  =\mathbb{E}\left(\varphi_k(U_0)\varphi_k(U_1)\varphi_j(U_m)\varphi_j(U_{m+1})\right)-\lambda_k\lambda_j= & \nonumber\\
=\lambda_k\lambda_j\sum\limits_{z=1}^\infty\lambda_z^{m-1}\left(\int_0^1\varphi_z(x)\varphi_k^2(x)dx\right)\left(\int_0^1\varphi_z(x)\varphi_j^2(x)dx\right). &
\end{eqnarray*}
Therefore, formula (\ref{covv}) becomes 
\begin{eqnarray*}\label{cov}
\underset{k\neq j}{\Sigma_{k,j}}&=&\sum\limits_{z=1}^s\lambda_z(\int_0^1\varphi_z(x)\varphi_k(x)\varphi_j(x)dx)^2+2\lambda_k\lambda_j(\int_0^1\varphi^2_k(x)\varphi_j^2(x)dx-1)\nonumber\\&&+2\lambda_k\lambda_j\sum\limits_{z=1}^s\sum\limits_{m=2}^\infty\lambda_z^{m-1}(\int_0^1\varphi_z(x)\varphi_k^2(x)dx)(\int_0^1\varphi_z(x)\varphi_j^2(x)dx)-\lambda_k\lambda_j.
\end{eqnarray*}

\subsection{Proof of Proposition \ref{components}}
Using the provided theorems for each of the copula families, we get:

{\it 1.  The sine copula},   $\varphi_1(x)=\sqrt{2}\cos \pi x$ ~and~ $\varphi_2(x)=\sqrt{2}\cos 2\pi x$. 
It follows that 
\begin{eqnarray}\label{cos} 
\sigma_1^2&=&1-\lambda_1^2 +\lambda_1(\int_0^1(\sqrt{2}\cos \pi x)^3dx)^2+\lambda_2(\int_0^1\sqrt{2}\cos 2\pi x(\sqrt{2}\cos \pi x)^2dx)^2+\nonumber\\&&+2\lambda_1^2(\int_0^1(\sqrt{2}\cos \pi x)^4dx-1)+2\lambda_1^2\sum\limits_{m=2}^\infty\lambda_1^{m-1}(\int_0^1(\sqrt{2}\cos \pi x)^3dx)^2\nonumber\\&&+2\lambda_1^2\sum\limits_{m=2}^\infty\lambda_2^{m-1}(\int_0^1\sqrt{2}\cos 2\pi x(\sqrt{2}\cos \pi x)^2dx)^2.\end{eqnarray}
Given that 
$\displaystyle
\int_0^1(\sqrt{2}\cos \pi x)^3dx)
= 0,
\int_0^1\sqrt{2}\cos 2\pi x (\sqrt{2}\cos \pi x)^2dx)
=\dfrac{1}{\sqrt{2}}$  and 
$$\int_0^14\cos^4 \pi xdx-1=\dfrac{1}{2}, \quad \text{ we obtain}\quad \sigma_1^2=1+\dfrac{\lambda_2}{2}+\lambda_1^2\dfrac{\lambda_2}{1-\lambda_2}.$$
$$\text{Similarly,}\quad \int_0^1(\sqrt{2}\cos 2\pi x)^3dx=0, \int_0^1\sqrt{2}\cos \pi x(\sqrt{2}\cos 2\pi x)^2dx)=0$$ and
$ \int_0^1(\sqrt{2}\cos 2\pi x)^4dx-1=\dfrac{1}{2}$
imply 
$
\sigma^2_2=1,
$
$\displaystyle \Sigma_{1,2}=\Sigma_{2,1}=\lambda_1\left(\dfrac{1}{2}-\lambda_2\right).$

{\it 2. The Sine-cosine copula.} We have
\[\varphi_1(x)=\cos 2\pi x;~~\varphi_2(x)=\cos 4\pi x;~~\varphi_3(x)=\sin 2\pi x \text{ and } \varphi_4(x)=\sin 4\pi x.\]

Via simple computations, we get
    \begin{eqnarray}\label{cal3}
    \int_0^1(\sqrt{2})^3\cos^3 2\pi x dx
    =0,\int_0^1(\sqrt{2})^3\cos^2 2\pi x \cos 4\pi x dx
    =\dfrac{1}{\sqrt{2}},&&~\nonumber\\\int_0^1 (\sqrt{2})^3\cos^2 2\pi x \sin 2\pi x dx
    =0,\int_0^12\sqrt{2}\cos^2 2\pi x \sin 4\pi x
    =0,&&\nonumber\\
     \text{ and } \int_0^14\cos^4 2\pi xdx=\dfrac{3}{2}. \quad \text{Using these integrals leads to}\nonumber
    \end{eqnarray}
$$\Sigma_{1,1}=\sigma^2_1=1-\lambda_1^2+\dfrac{\lambda_2}{2}+\lambda_1^2\lambda_1^2+\dfrac{\lambda_2}{1-\lambda_2}=1+\dfrac{\lambda_2}{2}+\lambda_1^2\dfrac{\lambda_2}{1-\lambda_2}.$$

    {\begin{eqnarray*}\bullet~ \Sigma_{2,2}=\sigma_2^2=1-\lambda_2^2 +\lambda_1(\int_0^12\sqrt{2}\cos 2\pi x \cos^2 4\pi xdx)^2+\\ \lambda_2(\int_0^12\sqrt{2}\cos^3 4\pi xdx)^2+\mu_1(\int_0^12\sqrt{2}\cos^2 4\pi x\sin 2\pi xdx)^2+\\ \mu_2(\int_0^12\sqrt{2}\cos^2 4\pi x\sin 4\pi xdx)^2+2\lambda_2^2(\int_0^14\cos^4 4\pi xdx-1)+\\ 2\lambda_2^2\sum\limits_{m=2}^\infty\lambda_1^{m-1}(\int_0^12\sqrt{2}\cos^2 4\pi x\cos 2\pi xdx)^2+\\2\lambda_2^2\left(\sum\limits_{m=2}^\infty\lambda_2^{m-1}(2\sqrt{2}\int_0^1\cos^3 4\pi xdx)^2\right)+    \\ +2\lambda_2^2\left(\sum\limits_{m=2}^\infty\mu_1^{m-1}(\int_0^12\sqrt{2}\sin 2\pi x\cos^2 4\pi xdx)^2\right)+\\2\lambda_2^2\sum\limits_{m=2}^\infty\mu_2^{m-1}(\int_0^12\sqrt{2}\sin 4\pi x\cos^2 4\pi xdx)^2.\end{eqnarray*}}
   Knowing $\displaystyle
    \int_0^1 2\sqrt{2}\cos^3 4\pi xdx=
     \int_0^1 2\sqrt{2}\cos 2\pi x \cos^2 4\pi x dx=0$, $\displaystyle
    \int_0^1 2\sqrt{2}\cos^2 4\pi x \sin 2\pi x dx=\int_0^1 2\sqrt{2}\cos^2 4\pi x \sin 4\pi x dx=
    0,$  $\text{and} \quad\displaystyle \int_0^12\cos^2 4\pi x dx-1=\dfrac{1}{2},$ we derive
    $ \Sigma_{2,2}=\sigma_2^2=1-\lambda_2^2+\lambda_2^2=1.$

{\scriptsize\begin{eqnarray*}
\bullet~ \Sigma_{3,3}=\sigma_3^2=1-\mu_1^2 +\lambda_1(\int_0^1\sqrt{2}\cos 2\pi x(\sqrt{2}\sin 2\pi x)^2dx)^2+\\
\lambda_2(\int_0^1\sqrt{2}\cos 4\pi x(\sqrt{2}\sin 2\pi x)^2dx)^2\nonumber+\mu_1(\int_0^1(\sqrt{2}\sin 2\pi x)^3dx)^2+\\
\mu_2(\int_0^1(\sqrt{2}\sin 2\pi x)^2\sqrt{2}\sin 4\pi xdx)^2+2\mu_1^2(\int_0^1(\sqrt{2}\sin 2\pi x)^4dx-1)+\\
2\mu_1^2\sum\limits_{m=2}^\infty\lambda_1^{m-1}(\int_0^1(\sqrt{2}\sin 2\pi x)^2\sqrt{2}\cos 2\pi xdx)^2+\\
2\mu_1^2\sum\limits_{m=2}^\infty\lambda_2^{m-1}(\int_0^1(\sqrt{2}\sin 2\pi x)^2\sqrt{2}\cos 4\pi xdx)^2+\\
2\mu_1^2\sum\limits_{m=2}^\infty\mu_1^{m-1}(\int_0^1(\sqrt{2}\sin 2\pi x)^3dx)^2+\\
2\mu_1^2\sum\limits_{m=2}^\infty\mu_2^{m-1}(\int_0^1\sqrt{2}\sin 4\pi x(\sqrt{2}\sin 2\pi x)^2dx)^2
\end{eqnarray*}}

   \begin{eqnarray*}\label{cal5} \text{and}\quad
    \int_0^1 2\sqrt{2}\sin^3 2\pi x dx
   = 0, \int_0^1 2\sqrt{2}\cos 2\pi x \sin^2 2\pi x dx
    =0, \\ \int_0^1 2\sqrt{2}\sin^2 2\pi x \cos 4\pi x dx
    =-\dfrac{1}{\sqrt{2}},\\ \int_0^1 2\sqrt{2}\sin^2 2\pi x \sin 4\pi x dx
    =0, \int_0^14\sin 2\pi xdx-1=\dfrac{1}{2}.
\end{eqnarray*}

    \begin{eqnarray*}\text{Thus},\quad \Sigma_{3,3}=\sigma_3^2&=&1-\mu_1^2+\dfrac{\lambda_2}{2}+\mu_1^2+\mu_1^2\dfrac{\lambda_2}{1-\lambda_2}=1+\dfrac{\lambda_2}{2}+\mu_1^2\dfrac{\lambda_2}{1-\lambda_2}.\end{eqnarray*}
   
{\scriptsize\begin{eqnarray*} \bullet\Sigma_{4,4}=\sigma_4^2=1-\mu_2^2 +\lambda_1(\int_0^1\sqrt{2}\cos 2\pi x\cdot(\sqrt{2}\sin 4\pi x)^2dx)^2+\\
\lambda_2(\int_0^1\sqrt{2}\cos 4\pi x\cdot(\sqrt{2}\sin 4\pi x)^2dx)^2+\mu_1(\int_0^1\sqrt{2}\sin 2\pi x\cdot(\sqrt{2}\sin 4\pi x)^2dx)^2+\\ 
\mu_2(\int_0^1(\sqrt{2}\sin 4\pi x)^3dx)^2+2\mu_2^2(\int_0^1(\sqrt{2}\sin 4\pi x)^4dx-1)+\\
2\mu_2^2\sum\limits_{m=2}^\infty\lambda_1^{m-1}(\int_0^1(\sqrt{2}\sin 4\pi x)^2\cdot\sqrt{2}\cos 2\pi xdx)^2+\\
2\mu_2^2\sum\limits_{m=2}^\infty\lambda_2^{m-1}(\int_0^1(\sqrt{2}\sin 4\pi x)^2\cdot\sqrt{2}\cos 4\pi xdx)^2+\\
2\mu_2^2\sum\limits_{m=2}^\infty\mu_1^{m-1}(\int_0^1\sqrt{2}\sin 2\pi x\cdot(\sqrt{2}\sin 4\pi x)^2dx)^2+\\
2\mu_2^2\sum\limits_{m=2}^\infty\mu_2^{m-1}(\int_0^1(\sqrt{2}\sin 4\pi x)^3dx)^2,\end{eqnarray*}}

\begin{eqnarray*}\text{and} \int_0^1\left(\sqrt{2}\sin 4\pi x\right)^3dx=0, \int_0^1\left(\sqrt{2}\cos 2\pi x\right)\left(\sqrt{2}\sin 4\pi x\right)^2dx=0, \\ \int_0^1\left(\sqrt{2}\sin 4\pi x\right)^2\left(\sqrt{2}\cos 4\pi x\right)dx=0, \\ \int_0^1\left(\sqrt{2}\sin 2\pi x\right)\left(\sqrt{2}\sin 4\pi x\right)^2dx=0, \int_0^1\left(\sqrt{2}\sin 4\pi x\right)^4dx-1=\dfrac{1}{2}.\end{eqnarray*}
 \[\text{Therefore,}\quad  \Sigma_{4,4}=\sigma_4^2=1-\mu_2^2+\mu_2^2=1.\]

    {\begin{eqnarray*}\bullet~\Sigma_{1,2}&=&\sum\limits_{z=1}^4\lambda_z(\int_0^1\varphi_z(x)\varphi_1(x)\varphi_2(x)dx)^2+\lambda_1\lambda_2\left(2(\int_0^1\varphi^2_1(x)\varphi_2^2(x)dx-1)\right.\nonumber\\&&\left.+2\sum\limits_{z=1}^4\sum\limits_{m=2}^\infty\lambda_z^{m-1}(\int_0^1\varphi_z(x)\varphi_1^2(x)dx)(\int_0^1\varphi_z(x)\varphi_2^2(x)dx)-1\right).
\end{eqnarray*}
Moreover, just as for the other cases, we have
    {\begin{eqnarray*}\label{cal7}
    \int_0^1\varphi^2_1(x)\varphi_2(x)dx=\dfrac{1}{\sqrt{2}}, \int_0^1\varphi^2_2(x)\varphi_1(x)dx= \int_0^1\varphi_3(x)\varphi_1(x)\varphi_2(x)dx=0,\\ \int_0^1\varphi_4(x)\varphi_1(x)\varphi_2(x)dx=
    \int_0^1\varphi^2_1(x)\varphi^2_2(x)dx-1= \int_0^1\varphi^2_1(x)\varphi_3(x)dx=0, \\
    \int_0^1\varphi^2_1(x)\varphi_4(x)dx=
   \int_0^1\varphi^2_2(x)\varphi_3(x)dx=0,
    \text{ and } \int_0^1\varphi^2_2(x)\varphi_4(x)dx=0.
    \end{eqnarray*}}

\[\text{Therefore}, \quad \Sigma_{1,2}=\Sigma_{2,1}=\dfrac{\lambda_1}{2}-\lambda_1\lambda_2=\lambda_1\left(\dfrac{1}{2}-\lambda_2\right).\]
    
 { \begin{eqnarray*}
\bullet~\Sigma_{2,3}&=&\sum\limits_{z=1}^4\lambda_z(\int_0^1\varphi_z(x)\varphi_3(x)\varphi_2(x)dx)^2+\lambda_2\mu_1\left(2(\int_0^1\varphi^2_3(x)\varphi_2^2(x)dx-1)\right.\nonumber\\&&\left.+2\sum\limits_{z=1}^4\sum\limits_{m=2}^\infty\lambda_z^{m-1}(\int_0^1\varphi_z(x)\varphi_3^2(x)dx)(\int_0^1\varphi_z(x)\varphi_2^2(x)dx)-1\right).\end{eqnarray*}}

    Moreover, 
\begin{eqnarray*}\label{cal8}
    \int_0^1\varphi^2_2(x)\varphi_3(x)dx=\int_0^1\varphi_4(x)\varphi_3(x)\varphi_2(x)dx=
    \int_0^1\varphi^2_3(x)\varphi^2_2(x)dx=0,\\
    \int_0^1\varphi^2_3(x)\varphi_1(x)dx=\int_0^1\varphi^2_3(x)\varphi_4(x)dx=0,
    \int_0^1\varphi^2_3(x)\varphi_2(x)dx=-\dfrac{1}{\sqrt{2}}.
\end{eqnarray*}

     \[\text{Theorefore,}\quad \Sigma_{2,3}=\Sigma_{3,2}=\dfrac{\mu_1}{2}-\lambda_2\mu_1=\mu_1\left(\dfrac{1}{2}-\lambda_2\right).\]   

 {\begin{eqnarray*}\bullet~\Sigma_{1,3}=\sum\limits_{z=1}^4\lambda_z(\int_0^1\varphi_z(x)\varphi_3(x)\varphi_1(x)dx)^2+\lambda_1\mu_1\left(2(\int_0^1\varphi^2_3(x)\varphi_1^2(x)dx-1)\right.\nonumber\\
\left.+2\sum\limits_{z=1}^4\sum\limits_{m=2}^\infty\lambda_z^{m-1}(\int_0^1\varphi_z(x)\varphi_3^2(x)dx)(\int_0^1\varphi_z(x)\varphi_1^2(x)dx)-1\right).\end{eqnarray*}}

     \begin{equation*}\text{Since}, \label{cal10}\int_0^1\varphi^2_3(x)\varphi^2_1(x)dx-1=-\dfrac{1}{2},
\int_0^1\varphi_4(x)\varphi_3(x)\varphi_1(x)dx=\dfrac{1}{\sqrt{2}},\end{equation*}
\[ \text{it follows that}\quad \Sigma_{1,3}=\Sigma_{3,1}=\dfrac{\mu_2}{2}-2\lambda_1\mu_1.\]

 {\begin{eqnarray*}\bullet~
\Sigma_{1,4}=\sum\limits_{z=1}^4\lambda_z(\int_0^1\varphi_z(x)\varphi_4(x)\varphi_1(x)dx)^2+\lambda_1\mu_2\left(2(\int_0^1\varphi^2_4(x)\varphi_1^2(x)dx-1)\right.\nonumber\\
\left.+2\sum\limits_{z=1}^4\sum\limits_{m=2}^\infty\lambda_z^{m-1}(\int_0^1\varphi_z(x)\varphi_4^2(x)dx)(\int_0^1\varphi_z(x)\varphi_1^2(x)dx)-1\right).\end{eqnarray*}}

{\begin{eqnarray*}\label{cal9} \text{ Moreover,}
    \int_0^1\varphi^2_4(x)\varphi^2_1(x)dx-1=
\int_0^1\varphi^2_4(x)\varphi_1(x)dx=0,\\
\int_0^1\varphi_4(x)\varphi_1(x)\varphi_3(x)dx=\dfrac{1}{\sqrt{2}}. \quad \text{Thus,}\quad \Sigma_{1,4}=\Sigma_{4,1}=\dfrac{\mu_1}{2}-\lambda_1\mu_2.
\end{eqnarray*}} 
     \begin{eqnarray*}
\bullet~ \Sigma_{2,4}&=&\sum\limits_{z=1}^4\lambda_z(\int_0^1\varphi_z(x)\varphi_4(x)\varphi_2(x)dx)^2+\lambda_2\mu_2\left(2(\int_0^1\varphi^2_4(x)\varphi_2^2(x)dx-1)\right.\nonumber\\&&\left.+2\sum\limits_{z=1}^4\sum\limits_{m=2}^\infty\lambda_z^{m-1}(\int_0^1\varphi_z(x)\varphi_4^2(x)dx)(\int_0^1\varphi_z(x)\varphi_2^2(x)dx)-1\right).
\end{eqnarray*}
    \begin{eqnarray*}\label{cal11}\text{So,}
    \int_0^1\varphi^2_4(x)\varphi^2_2(x)dx=-\dfrac{1}{2},
\int_0^1\varphi^2_4(x)\varphi_2(x)dx=0,  \Sigma_{2,4}=\Sigma_{4,2}=-2\lambda_2\mu_2. \end{eqnarray*}
  
   {\begin{eqnarray*}
\bullet~ \Sigma_{3,4}&=&\sum\limits_{z=1}^4\lambda_z(\int_0^1\varphi_z(x)\varphi_4(x)\varphi_3(x)dx)^2+\mu_1\mu_2\left(2(\int_0^1\varphi^2_4(x)\varphi_3^2(x)dx-1)\right.\nonumber\\&&\left.+2\sum\limits_{z=1}^4\sum\limits_{m=2}^\infty\lambda_z^{m-1}(\int_0^1\varphi_z(x)\varphi_4^2(x)dx)(\int_0^1\varphi_z(x)\varphi_3^2(x)dx)-1\right).\end{eqnarray*}}
    \begin{eqnarray*}
    \text{But,}\quad \int_0^1\varphi^2_4(x)\varphi^2_3(x)dx-1=
    \int_0^1\varphi^2_4(x)\varphi_3(x)dx=0,\nonumber\\
\int_0^1\varphi^2_4(x)\varphi_2(x)dx=0,
\int_0^1\varphi_4(x)\varphi_1(x)\varphi_3(x)dx=\dfrac{1}{\sqrt{2}}.\end{eqnarray*}
     \[\text{Thus,}\quad \Sigma_{2,4}=\Sigma_{4,2}=\dfrac{\lambda_1}{2}-\mu_1\mu_2.\]

{\it For the Legendre copula:} 
{\begin{eqnarray*}
\int_0^1\varphi^3_1(x)dx)=\int_0^1(\sqrt{3}(2x-1))^3dx=0,\nonumber\\
\int_0^1\varphi_2(x)\varphi_1^2(x)dx)=\int_0^1\sqrt{5}(6x^2-6x+1) (\sqrt{3}(2x-1))^2dx=\dfrac{2}{\sqrt{5}},\nonumber\\
\int_0^1\varphi^4_1(x)dx-1=\int_0^1(\sqrt{3}(2x-1))^4dx-1=\dfrac{4}{5}.
\end{eqnarray*}}
\[ \text{Thus,}\quad \sigma_1^2=1-\lambda_1^2 +\dfrac{4}{5}\lambda_2+\dfrac{8}{5}\lambda_1^2+\dfrac{8}{5}\lambda_1^2\sum\limits_{m=2}^\infty\lambda_2^{m-1}=1+\dfrac{4}{5}\lambda_2+\lambda_1^2\dfrac{3+5\lambda_2}{5(1-\lambda_2)} .\]

\begin{eqnarray*} \text{Similarly,}\quad
\int_0^1\varphi^3_2(x)dx)=\int_0^1(\sqrt{5}(6x^2-6x+1))^3dx=\dfrac{2\sqrt{5}}{7},\nonumber\\
\int_0^1\varphi_1(x)\varphi_2^2(x)dx)=\int_0^1 5\sqrt{3}(6x^2-6x+1)^2 (2x-1)dx=0,\nonumber\\
\int_0^1\varphi^4_2(x)dx-1=\int_0^1(\sqrt{5}(6x^2-6x+1))^4dx-1=\dfrac{8}{7},\quad \text{implies}
\end{eqnarray*}}

\begin{eqnarray*}
\sigma_2^2&=&1-\lambda_2^2 +\dfrac{20}{49}\lambda_2+\dfrac{16}{7}\lambda_2^2+\dfrac{40}{49}\lambda_2^2\sum\limits_{m=2}^\infty\lambda_2^{m-1}=1+\dfrac{20}{49}\lambda_2+\lambda_2^2\dfrac{63-23\lambda_2}{49(1-\lambda_2)} .\end{eqnarray*}
The correlation in this case is obtained as
 \begin{eqnarray*}
\Sigma_{1,2}&=&\sum\limits_{z=1}^2\lambda_z(\int_0^1\varphi_z(x)\varphi_1(x)\varphi_2(x)dx)^2+\lambda_1\lambda_2\left(2(\int_0^1\varphi^2_1(x)\varphi_2^2(x)dx-1)\right.\nonumber\\&&\left.+2\sum\limits_{z=1}^2\sum\limits_{m=2}^\infty\lambda_z^{m-1}(\int_0^1\varphi_z(x)\varphi_1^2(x)dx)(\int_0^1\varphi_z(x)\varphi_2^2(x)dx)-1\right).
\end{eqnarray*}

{\begin{eqnarray*}\quad\text{But,}
\int_0^1\varphi_2(x)\varphi_1^2(x)dx)=\dfrac{2}{\sqrt{5}},
\int_0^1\varphi_1^2(x)\varphi^2_2(x)dx-1=\dfrac{4}{7}.
\end{eqnarray*}}

\begin{eqnarray*} \text{So,}\quad
\Sigma_{1,2}&=&\dfrac{4}{5}\lambda_1+\lambda_1\lambda_2\left(\dfrac{8}{7}-1+\dfrac{8}{7}\dfrac{\lambda_2}{(1-\lambda_2)}\right)=\dfrac{4}{5}\lambda_1+\lambda_1\lambda_2\dfrac{1+7\lambda_2}{7(1-\lambda_2)}.\end{eqnarray*}

\subsection{Proof of Theorem \ref{MLE1}}

\underline{Verification of conditions}
\begin{enumerate}
\item Considering the parameters space,  the density of each of the copulas is $c(u,v)\le M<2$ inside the parameter space. Therefore, according to Longla et al. (2022c), the Markov chains generated by such copulas are $\psi-$mixing. Since $\psi$-mixing implies ergodicity, the first condition is satisfied.
\item The second condition is satisfied as well since $c(u,v)>0$  on $[0,1]\times [0,1]$ for any set of parameters not on the boundary of the parameter set.
\item Since $c(u,v)>0$ and is continuous on $[0,1]^2$, all partial derivatives mentioned in the third condition exist and are continuous on $[0,1]^2$ when $\Lambda$ is not on the boundary of the set.
\item Since there is a constant $M<c(u,v)<2$ on $[0,1]^2$ when $\Lambda$ is an interior point, then
 \[ \mathbb{E}_\Lambda\left[\left(\frac{\partial}{\partial\lambda_i}\ln c(u,v)\right)^2\right]<\infty.\]
Via simple calculations, we establish that when the parameters are strictly inside the region, we have for some $K$ not depending on $\Lambda$,  
\[\lvert\frac{\partial^3}{\partial\lambda_i\lambda_j\lambda_k}\ln c(u,v)\lvert\leq K \quad \text{and}\quad
 \mathbb{E}_\Lambda\left[\underset{\Lambda\in A'}{\sup}\lvert\frac{\partial^3}{\partial\lambda_i\lambda_j\lambda_k}\ln c(u,v)\lvert\right]<\infty.  \]

\item Since entries of the Fisher information matrix $\Sigma(\Lambda)$ cannot be computed in closed form, an estimation procedure is used in Matlab to approximate it. It appears to have a strictly positive determinant for each of the cases, implying that it is inversible.
Considering the extra regularity condition, $n\Sigma$ is aproximated as an average by $I_n(\Lambda)$ because $I_n/n$ converges in probability to $\Sigma$ thanks to ergodicity (see Sun et al (2020)). Here,
     \[ I_n(\Lambda)=-\left(\frac{\partial^2}{\partial\lambda_k\lambda_j}~\ell(\Lambda)\right)_{1\leq k,j\leq s}. \]
\end{enumerate}
This completes the proof of the following statement. $\square$

\end{document}